\documentclass[12pt]{amsart}
\usepackage{amssymb}

\begin{document}
\baselineskip 18pt
\hfuzz=6pt

\newtheorem{theorem}{Theorem}[section]
\newtheorem{prop}[theorem]{Proposition}
\newtheorem{lemma}[theorem]{Lemma}
\newtheorem{definition}[theorem]{Definition}
\newtheorem{cor}[theorem]{Corollary}
\newtheorem{example}[theorem]{Example}
\newcommand{\ra}{\rightarrow}
\renewcommand{\theequation}
{\thesection.\arabic{equation}}
\newcommand{\ccc}{{\cal C}}
\newcommand{\one}{1\hspace{-4.5pt}1}

\def\BRR  { ${\rm  BMO}({\mathbb R} ) $}
\def\BRRn  { ${\rm  BMO}({\mathbb R}^n )$}
\def\BR  {{${\rm  BMO}$}}
\def\br {\rm BMO}
\def\HLR  { $H^{1}_{L}({\mathbb R} ) $}
\def\BLR  { ${\rm  BMO}_{L}({\mathbb R} ) $}
\def\BLO  { ${\rm  BMO}_{L}({\Omega} ) $}
\def\BLRD  { ${\rm  BMO}_{{L}^{\ast}}  $}
\def\HL  { $H^{1}_{L} $}
\def\HLD  { $H^{1}_{L} \cap L^2  $}
\def\HR  {{$H^1({\mathbb R} ) $}}
\def\BL  {${\rm  BMO}_{L}$}
\def\RR{\mathbb R}

\title[Comparison  of   the classical {BMO} with the new {BMO}  spaces ]
{Comparison  of   the classical {BMO} with the {BMO}  spaces
      associated  with operators  and applications}
\thanks{{\it 2000 Mathematical Subject Classification:} 42B35, 42B25, 47B38.}
\thanks{{\it Key words and phrases:} {\BR} space, Hardy space, 
Dirichlet and Neumann
Laplacian, semigroup, Gaussian bounds,
fractional powers, purely imaginary powers, spectral multiplier.}
\thanks{D.G. Deng and L.X. Yan are  supported
by NNSF of China  and  the Foundation of  Advanced Research
Center, Zhongshan University.
X.T. Duong and L.X. Yan are supported by a grant from
Australia Research Council. Also A. Sikora participated in the collaboration
with the support of the ARC}
\date{}

\author{{ {Donggao Deng,\   Xuan Thinh Duong,\
Adam Sikora  \ and \  Lixin Yan}}}

\begin{abstract}
Let $L$ be a generator of  a semigroup satisfying the Gaussian upper bounds.
In this paper, we study further a new {\BL}   space associated with $L$
which was
introduced recently by Duong and Yan. We discuss applications
of the new {\BL} spaces in the theory of singular integration such as
{\BL} estimates and interpolation results for fractional powers,
purely imaginary powers  and spectral multipliers of self adjoint operators.
We also demonstrate that the space {\BL} might coincide with or might be
essentially
different from the classical {\BR} space.

\end{abstract}


\medskip


\maketitle

\bigskip

\section{Introduction}
\setcounter{equation}{0}
The classical space of functions of bounded mean oscillation  ({\BR})
plays a crucial role in  modern harmonic analysis.
See for examples  \cite{FS}, \cite{JN}, \cite{St1} and \cite{St2}.
In the case of the
Euclidean space $\mathbb R^n$, a function $f$ is said to in{\BRRn}
if
\begin{eqnarray}
\label{e2.3a}
      \|f\|_{\br({\mathbb R}^n)}=\sup_Q\frac{1}{|Q|}\int_Q|
f(x)-f_Q|dx <\infty,
\end{eqnarray}
where $f_Q$ denotes the average value of $f$  on the cube $Q$ and
the supremum is taken over all cubes $Q$  in ${\mathbb R^n} $.

An important  application of the theory of  {\BR}
spaces is  the following
interpolation result.

\begin{prop}\label{prop1}
If $T$ is a bounded sublinear
operator from $L^2({\mathbb R}^n)$ to $L^2({\mathbb R}^n)$, and $T$ is
bounded from $L^{\infty}({\mathbb R}^n)$ to {\BRRn},
then
$T$ is bounded from $L^p({\mathbb R}^n)$ to $L^p({\mathbb R}^n)$ for
all $2 < p < \infty$.
\end{prop}

It is well known that Calder\'on-Zygmund   operators (such as
the Hilbert transform on the real line, the Riesz transforms on
${\mathbb R}^n$, or the purely imaginary powers of the Laplacian
on ${\mathbb R}^n$) do not map the space $L^{\infty}$ into
$L^{\infty}$, but the standard conditions on their kernels ensure that
they map $L^{\infty}$ into the {\BR} space boundedly, hence we
can apply Proposition~\ref{prop1} to obtain $L^p$ boundedness of
these operators
for $p>2$.
In this sense, the BMO space  is a natural
substitute of the space $L^{\infty}$ in the  theory   of Calder\'on-Zygmund
    singular integrals.

In this paper we are motivated by study of singular integral operators 
corresponding to spectral multiplier of an operator $L$ 
which generates a semigroup with appropriate  kernel bounds, see
\cite{DY1}. Such multipliers do not always map $L^\infty$ or appropriate
$L^p$ spaces into the classical  {\BR} space, see Example~\ref{theo45} below.
Hence  the classical {\BR} space is not  necessarily a
suitable space to study such singular integrals. 
To study these rough operators,
we introduced a new {\BL}\ space associated with an operator
$L$.

To explain our approach to {\BL}\ space associated with an operator
let us recall that  the space of {\BR} functions can be characterized by the Carleson measure
estimate as follows:

\begin{prop}\label{prop2}
    A function $f$ is in  {\BR}   if and
only if $f$ satisfies
$\int_{{\mathbb R}^n} \frac{|f(x)|}{1+|x|^{n+1}}dx<\infty$, and $\mu_f(x,t)=
|t\frac{\partial}{ \partial t} e^{-t\sqrt{{\Delta}}} f(x)|^2\frac{dxdt}{t}$
is a  Carleson measure.
\end{prop}

One can see from the characterization in Proposition~\ref{prop2} that the
{\BR} space is associated with the Laplace operator on ${\mathbb R}^n$
and it seems to be natural idea to replace the Laplace operator ${\Delta}$
by more general operators operator $L$, see also
\cite{FS} and \cite{St2}.
In this paper we use equivalent approach, see Definition~\ref{defi23} below.
 In this definition  the {\BL}\ space associated with $L$ is defined
by using the function
$e^{-t_QL}f$ to replace the average
$f_Q$ in  definition \ref{e2.3a} of {\BR}  where the value $t_Q$ is scaled
to the  length of the sides of Q.
In this paper we discuss various examples which shows that
Definition~\ref{defi23} is an effective tool in study of singular integrals 
operators associated with the operator $L$. 
We refer the reader to
\cite{ADM}, \cite{CW} and \cite{DGMTZ} for other ideas related to generalization
of the {\BR} space and BMO spaces associated with an operator $L$.

Many important features of the classical BMO space are retained
by the new {\BL} spaces such as the John-Nirenberg inequality and
duality between the Hardy space and the {\BL} space. See \cite{DY1}
and \cite{DY2}.
One of these important features is that the interpolation property in
Proposition~\ref{prop1} is still valid if the classical space {\BR}
is replaced by
the  {\BL} space associated with an operator $L$. Indeed, the following result
 is proved in \cite{DY1} (Theorem 6.1).

\begin{prop}\label{prop3}
    Let ${\mathcal X}$ be a space of homogeneous type.
If $T$ is a
bounded sublinear  operator from $L^2({\mathcal X})$ to
$L^2({\mathcal X})$, and $T$ is bounded
from $L^{\infty}({\mathcal X})$ into ${\rm BMO}_L({\mathcal X})$, then
$T$ is bounded from $L^p({\mathcal X})$ to $L^p({\mathcal X})$ for
all $2 < p < \infty$.
\end{prop}

A natural question arising from Proposition~\ref{prop3} is to compare the
classical BMO space and the  BMO$_L$ space  associated with an
operator $L$. 
In Sections 3 and 4 we study this question
systematically and  we show that depending on the choice of the
operator $L$, all the following cases are possible 

Case 1: BMO $\cong$ BMO$_L$;

Case 2: BMO $\subseteq$ BMO$_L$ and BMO $\neq$ BMO$_L$;

Case 3: BMO$_L$ $\subseteq$ BMO and BMO$_L$ $\neq$ BMO;

Case 4: ${\rm BMO}\not\subseteq {\rm BMO}_L$
and ${\rm BMO}_L\not\subseteq {\rm BMO}.$

    \smallskip

For other results related to Cases 1 and 2 we refer readers to
 Proposition 2.5 of \cite{DY1}, Section 6.2 of \cite{DY2}
and Proposition~3.1 of \cite{Ma}.
In Section 5 we show that if $f\in L^{n/\alpha}({\mathbb R}^n)$ 
and $L^{-\alpha}f <\infty$
almost everywhere then $L^{-\alpha}f\in  $BMO$_L$.
We construct an example of a function $f\in L^p({\mathbb R})$ and an
operator $L$ such that $L^{-\frac{1}{2p}}f \in$ BMO$_L$
but $L^{-\frac{1}{2p}} f \notin$ BMO. This shows that the new BMO$_L$ space
does make a difference in estimates of singular integrals.
Finally in Sections 6 and 7, we  obtain sharp estimates of the $L^\infty$ to  {\BL}
 norm of the  purely
imaginary powers
$L^{is}$ of a self adjoint operator $L$. We also obtain the BMO type estimates for 
spectral multipliers of a self adjoint operator $L$
and for  maximal operators $\sup_{t>0}|F(tL)|$ corresponding to $L$
and appropriate functions $F$.
$L^p$ boundedness of these operators, $2<p<\infty$, then follows from
Proposition 1.3.

\section{Preliminaries}
\setcounter{equation}{0}

\subsection{BMO spaces on the  half spaces.}

      Let us begin by  recalling the
definitions  of
various BMO spaces on the usual upper-half space  in ${\mathbb R}^n$.
For any subset $A\subset {\mathbb R}^n$ and a function 
$f \colon  {\mathbb R}^n   \to  {\mathbb C}$ by $ f|_A$ we denote 
the restriction of $f$ to the set $A$. Next we set
$$
{\mathbb R}^n_+=\Big\{(x', x_n)\in {\mathbb R}^n:
      \ x'=(x_1, \cdots, x_{n-1})\in {\mathbb R}^{n-1}, \ x_n>0\Big\}.
$$
\begin{definition}
\label{defi21}
A function $f$ on ${\mathbb R}^n_+$ is said to be in ${\rm
BMO}_r({\mathbb R}^n_+)$ if there exists $F\in {\rm
BMO}({\mathbb R}^n)$ such that $F|_{{\mathbb R}^n_+}=f$. If $f\in  {\rm
BMO}_r({\mathbb R}^n_+)$, we set
$$
\|f\|_{{\rm
BMO}_r({\mathbb R}^n_+)} =\inf \left\{\|F\|_{{\rm BMO}({\mathbb R}^n) }
   \colon     F|_{{\mathbb R}^n_+}=f   \right\}.
$$

      \medskip

A function $f$ on ${\mathbb R}^n_+$ belongs to
${\rm BMO}_z({\mathbb R}^n_+)$
      if
the function $F$
defined by
\begin{eqnarray}
\label{e2.6}
F(x)=
\left\{\begin{array}{ll}
f(x) \ \ & {\rm if} \ x\in{\mathbb R}^n_+;\\[5pt]
0 \ & {\rm if}\ x\not\in {\mathbb R}^n_+
\end{array}
\right.
\end{eqnarray}
belongs to ${\rm BMO}({\mathbb R}^n)$. If $f\in
{\rm BMO}_z({\mathbb R}^n_+)$, we set
       $\|f\|_{{\rm BMO}_z({\mathbb R}^n_+) } =
\|F\|_{{\rm BMO}_z({\mathbb R}^n) }.$
\end{definition}
Compare Section~4.5.1, page 221 of  \cite{Tr} and Section~5.4 of
\cite{ART}.
In order to analyze the spaces ${\rm BMO}_r({\mathbb R}^n_+)$ and ${\rm
BMO}_z({\mathbb R}^n_+)$, let us   introduce   the following notations,
see \cite{CKS}. For any $x=(x', x_n)\in {\mathbb R}^n $, we set
$
{\widetilde x}= (x', -x_n).
$
If $f$ is any function defined on ${\mathbb R}^n_+$, its even extension
$f_e$ is  defined on ${\mathbb R}^n$ by
\begin{eqnarray*}
f_{e}(x)=
\left\{\begin{array}{ll}
f(x)  \ \ & {\rm if} \ x\in {\mathbb R}^n_+;\\[5pt]
f({\widetilde x})  \ & {\rm if}\ x\in  {\mathbb R}^n_-,
\end{array}
\right.
\end{eqnarray*}
and its odd  extension    $f_{o}$
is defined by
\begin{eqnarray*}
f_{o}(x)=
\left\{\begin{array}{ll}
f(x)  \ \ & {\rm if} \ x\in {\mathbb R}^n_+;\\[5pt]
-f({\widetilde x})  \ & {\rm if}\ x\in {\mathbb R}^n_-,
\end{array}
\right.
\end{eqnarray*}
where
$$
{\mathbb R}^n_-=\Big\{(x', x_n)\in {\mathbb R}^n:  \ x'=(x_1, \cdots, x_{n-1})
\in {\mathbb R}^{n-1}, x_n<0\Big\}.
$$
For any function $f \in L^1_{\rm loc}({\mathbb R}^n_+),$
we define
\begin{eqnarray*}
\|f\|_{{\rm BMO}_e({\mathbb R}^n_+)} =\|f_e\|_{{\rm BMO}({\mathbb R}^n)}
\quad \mbox{and} \quad
\|f\|_{{\rm BMO}_o({\mathbb R}^n_+)} =\|f_o\|_{{\rm BMO}({\mathbb R}^n)}
\end{eqnarray*}
and we denote by ${{\rm BMO}_e({\mathbb R}^n_+)}$ and ${{\rm BMO}_o({\mathbb
R}^n_+)}$  the corresponding Banach spaces.
      We will see that ${\rm BMO}_e({\mathbb R}^n_+)$
is suitable for the analysis of the Neumann Laplacian on ${\mathbb R}^n_+$
whereas  ${\rm BMO}_o({\mathbb R}^n_+)$
is suitable for the study of the Dirichlet Laplacian on ${\mathbb R}^n_+.$
See Proposition  3.2 below.

In what follows,
$Q=Q[x_Q, l_Q]$   denotes  a cube  of ${\mathbb R}^n$
centered at $x_Q$ and of the side length
$l_Q.$
      Given any  cube $Q$, we denote the reflection of $Q$ across
$\partial {\mathbb R}^n_+$ by
      \begin{eqnarray}
       \widetilde{Q}=\Big\{(x', x_n)\in {\mathbb R}^n, \ \ (x', 
-x_n)\in Q\Big\}.
\label{e2.7}
\end{eqnarray}
Let $Q_+=Q\cap {\mathbb R}^n_+ $ and
      $Q_-=Q\cap {\overline{{\mathbb R}^n_-}}$ where $
{\overline{{\mathbb R}^n_-}}
=\big\{(x', x_n)\in {\mathbb R}^n:  \ x'=(x_1, \cdots, x_{n-1})
\in {\mathbb R}^{n-1}, x_n\leq 0\big\}.
$
If both $Q_-$ and $Q_+$ are not empty, we then   define
      \begin{eqnarray}
      \left\{
\begin{array}{ll}
{\widehat Q}_-=\{(x', x_n): x'\in Q\cap {\mathbb R}^{n-1},  \  -l_Q <x_n\leq
0\},\\[5pt]
{\widehat Q}_+=\{(x', x_n):\ x'\in Q\cap{\mathbb R}^{n-1}, \
0<x_n\leq   l_Q \}.
\label{e2.8}
\end{array}
\right.
\end{eqnarray}
Obviously,   we have the following properties:
   {\rm (i)}\ $Q_-\subseteq {\widehat Q}_-$,
$Q_+\subseteq {\widehat Q}_+$
and thus $Q\subseteq ({\widehat Q}_-\cup {\widehat Q}_+)$; {\rm (ii)}\
      $|Q|=|{\widehat Q}_-|=|{\widehat Q}_+|$.
      These will be often used in the sequel.

\subsection{Dirichlet and Neumann Laplacians.}

By ${\Delta}_{n, N_+}$ (and ${\Delta}_{n, N_-}$) we denote the
Neumann Laplacian
on ${\mathbb R}^n_+$ (and on ${\mathbb R}^n_-$ respectively).
Similarly by   ${\Delta}_{n, D_+}$ (and ${\Delta}_{n, D_-}$) we denote the
Dirichlet Laplacian
on ${\mathbb R}^n_+$ (and on ${\mathbb R}^n_-$ respectively).

The Dirichlet and Neumann Laplacian are positive definite self-adjoint
    operators. By the spectral theorem one can define
the semigroups generated by these operators
$\{\exp{(-t{\Delta}_{n, D_+})}\colon \,{t\geq 0} \}$ and
$\{\exp{(-t{\Delta}_{n, N_+}})\colon \, {t\geq 0}\}$.
By $p_{t,\ \!{\Delta}_{n, D_+}}(x,y)$ and  $p_{t, \ \!{\Delta}_{n,
N_+}}(x,y)$ we denote
the heat kernels corresponding to the semigroups generated by
${{\Delta}_{n, D_+}}$ and  ${{\Delta}_{n, N_+}}$ respectively.

For $n=1$ by the reflection method (see for example \cite[(6) p.
57]{S}) we obtain
\begin{eqnarray*}
p_{t, \ \! {\Delta}_{1, D_+}}(x,y)=\frac{1}{ (4\pi t)^{{1}/{2}}}
\Big(e^{-\frac{|x_1-y_1|^2}{ 4t}} -
e^{-\frac{|x_1+y_1|^2}{ 4t}}\Big).
\end{eqnarray*}
Then for $n \ge 2$
\begin{eqnarray}
p_{t,\ \! {\Delta}_{n, D_+}}(x,y)&=&\Big(p_{t, \ \! {\Delta}_{1,
D_+}}(x_n,y_n)\Big)
\Big(p_{t, \ \!{\Delta}_{n-1}}(x',y')\Big)\nonumber\\
&=&
\frac{1}{ (4\pi t)^{\frac{n}{2}}}  e^{-\frac{|x'-y'|^2}{ 4t}}
\Big(e^{-\frac{|x_n-y_n|^2}{ 4t}} -
e^{-\frac{|x_n+y_n|^2}{ 4t}}\Big),
\label{e2dkern}
\end{eqnarray}
where $p_{t, \ \! {\Delta}_{n-1}}(x,y)$ is the heat kernel corresponding to the
standard Laplace operator acting on ${\mathbb R}^{n-1}$.
Applying  the reflection method also to the Neumann Laplacian
we obtain (see \cite[(7) p. 57]{S})
\begin{eqnarray}
p_{t, \ \!{\Delta}_{n, N_+}}(x,y)&=&\Big(p_{t, \ \!{\Delta}_{1,
N_+}}(x_n,y_n)\Big)
\Big(p_{t,  \  \!{\Delta}_{n-1}}(x',y')\Big)\nonumber\\
&=&
\frac{1}{ (4\pi t)^{\frac{n}{2}}}  e^{-\frac{|x'-y'|^2}{ 4t}}
\Big(e^{-\frac{|x_n-y_n|^2}{ 4t}} +
e^{-\frac{|x_n+y_n|^2}{4t}}\Big).
\label{e2nkern}
\end{eqnarray}
In the sequel we skip the index $n$ and we denote the Dirichlet and Neumann
Laplacian  by  ${\Delta}_{D_+}$ and ${\Delta}_{N_+}$.
Note that by (\ref{e2dkern})
\begin{eqnarray}
\label{e2.12}
\exp(-t{\Delta}_{D_+}) f(x)&=&\int_{{\mathbb R}^n_+}
p_{t,  \  \!{\Delta}_{D_+}}(x,y)f(y)dy\nonumber\\&=&
\frac{1}{ (4\pi t)^{n/ 2}}\int_{{\mathbb R}^n}
e^{-\frac{|x-y|^2}{ 4t}} f_o(y)dy\nonumber\\
&=&\exp(-t{\Delta}) f_o(x)
\end{eqnarray}
for $x\in{\mathbb R}^n_+$  and all $t>0$.
Similarly
\begin{eqnarray}
\label{e2.14}
    \exp(-t{{\Delta}_{N_+}})f(x)&=&\int_{{\mathbb R}^n_+}
p_{t,  \  \!{\Delta}_{N_+}}(x,y)f(y)dy\nonumber\\
&=&\frac{1}{ (4\pi t)^{{n}/{2}}}\int_{{\mathbb R}^n}
e^{-\frac{|x-y|^2}{ 4t}} f_e(y)dy\nonumber\\
&=&\exp(-t{\Delta}) f_e(x)
\end{eqnarray}
for $x\in{\mathbb R}^n_+$  and all $t>0$.

Next for any  function $f$ on ${\mathbb R}^n$, we set
      $$
      f_-=f|_{{\mathbb R}^n_-}\quad {\rm and}\quad
      f_+=f|_{{\mathbb R}^n_+}.
$$
Now let ${\Delta}_N$ be the uniquely determined  unbounded operator
acting on   $L^2({\mathbb R}^{n})$  such that
\begin{equation}\label{abla}
({\Delta}_N f)_+={\Delta}_{N_+}f_+
\quad \mbox{and} \quad ({\Delta}_N f)_-={\Delta}_{N_-}f_-
\end{equation}
for all $f \colon {\mathbb R}^n \to  {\mathbb R}$ such that
$f_+ \in W^{1,2}({\mathbb R}^n_+)$ and $f_- \in W^{1,2}({\mathbb R}^n_-)$.
Then, ${\Delta}_N$ is a positive definite self-adjoint operator.
By (\ref{abla})
\begin{equation}\label{e2.19n}
(\exp(-t{\Delta}_N) f)_+=\exp(-t{\Delta}_{N_+})f_+
\quad \mbox{and} \quad
(\exp(-t{\Delta}_N) f)_-=\exp(-t{\Delta}_{N_-})f_-.
\end{equation}
Let $p_{t,  \ \!{\Delta}_{N}}(x,y)$ be the  heat kernel of
$\exp(-t{\Delta}_N)$. By
(\ref{e2.19n}) and (\ref{e2nkern}) we obtain
\begin{eqnarray}
p_{t,  \  \!{\Delta}_{N}}(x,y)
=\frac{1}{ (4\pi t)^{n/2}}  e^{-\frac{|x'-y'|^2}{ 4t}}
\Big(e^{-\frac{|x_n-y_n|^2}{ 4t}}+
e^{-\frac{|x_n+y_n|^2}{ 4t}}\Big) H(x_ny_n),
\label{e2.19}
\end{eqnarray}
where $H\colon \, \RR \to \{0,1\}$ is the Heaviside function  given by
\begin{eqnarray}
H(t)=
\left\{\begin{array}{ll}
0 \ \ & {\rm if} \ t< 0;\\\\
1 \ & {\rm if}\ t\geq  0.
\end{array}
\right.
\label{e2.20}
\end{eqnarray}
Similarly we define the Dirichlet Laplacian on ${\mathbb R}^{n}$ by the formula
\begin{equation}\label{diri}
({\Delta}_D f)_+={\Delta}_{D_+}f_+
\quad \mbox{and} \quad ({\Delta}_D f)_-={\Delta}_{D_-}f_-
\end{equation}
for all $f \colon {\mathbb R}^n \to  {\mathbb R}$ such that
$f_+ \in W_0^{1,2}({\mathbb R}^n_+)$ and $f_- \in W_0^{1,2}({\mathbb R}^n_-)$.
Then, ${\Delta}_D$ is a positive definite self-adjoint operator.
By (\ref{diri})
\begin{equation}\label{e2.19d}
(\exp(-t{\Delta}_D) f)_+=\exp(-t{\Delta}_{D_+})f_+
\quad \mbox{and} \quad
(\exp(-t{\Delta}_D) f)_-=\exp(-t{\Delta}_{D_-})f_-.
\end{equation}
Hence by  (\ref{e2dkern})  the kernel
$p_{t,  \  \! {\Delta}_{D}}(x,y)$  of the operator $\exp(-t{\Delta}_D)$ is
given by 
\begin{eqnarray}
p_{t,  \  \!{\Delta}_{D}}(x,y)=\frac{1}{ (4\pi t)^{n/ 2}}
e^{-\frac{|x'-y'|^2}{ 4t}}
\Big(e^{-\frac{|x_n-y_n|^2}{4t}}-
e^{-\frac{|x_n+y_n|^2}{4t}}\Big) H(x_ny_n).
\label{e2.22}
\end{eqnarray}

Finally we define the Dirichlet-Neumann Laplacian by the formula
\begin{equation}\label{dn}
({\Delta}_{DN} f)_+={\Delta}_{N_+}f_+
\quad \mbox{and} \quad ({\Delta}_{DN} f)_-={\Delta}_{D_-}f_-
\end{equation}
for all $f \colon {\mathbb R}^n \to  {\mathbb R}$ such that
$f_+ \in W^{1,2}({\mathbb R}^n_+)$ and $f_- \in W_0^{1,2}({\mathbb R}^n_-)$.
By (\ref{dn})
\begin{equation}\label{e2.19dn}
(\exp(-t{\Delta}_{DN}) f)_+=\exp(-t{\Delta}_{N_+})f_+
\,\, \mbox{and} \, \,
(\exp(-t{\Delta}_{DN}) f)_-=\exp(-t{\Delta}_{D_-})f_-.
\end{equation}
Hence by (\ref{e2dkern}) and (\ref{e2nkern}), the kernel
$p_{t,  \  \! {\Delta}_{DN}}(x,y)$  of $\exp(-t{\Delta}_{DN})$ is given by
\begin{equation}
p_{t,  \  \!{\Delta}_{DN}}(x,y)=\frac{1}{ (4\pi t)^{n/ 2}}
e^{-\frac{|x'-y'|^2}{ 4t}}
\Big(e^{-\frac{|x_n-y_n|^2}{ 4t}}+(2H(x_n)-1)
e^{-\frac{|x_n+y_n|^2}{ 4t}}\Big) H(x_ny_n).
\label{e2.22nd}
\end{equation}

      \medskip

Let us note that

      \medskip

       ($\alpha$)   All the  operators
${\Delta}$, ${\Delta}_{N_+}$,  ${\Delta}_{D_+}$,
      ${\Delta}_{N_-}$,  ${\Delta}_{D_-}$
and ${\Delta}_{D}$,  ${\Delta}_{N}$,
${\Delta}_{DN}$
are  self-adjoint and they generate bounded analytic positive semigroups
   acting on all $L^p$ spaces for $1\leq
p\leq \infty$;

      \medskip

($\beta$) Suppose that  $p_{t, L}(x,y)$ is  the kernel corresponding to
the semigroup generated by $L$ and that $L$ is one of the operators listed in
($\alpha$). Then the kernel $p_{t, L}(x,y)$ satisfies Gaussian bounds,
that is
\begin{eqnarray}
|p_{t, L}(x,y)|\leq \frac{C}{t^{n/2}} e^{-c\frac{|x-y|^2}{t}}
\label{e2.1}
\end{eqnarray}
for all $x,y \in \Omega$, where $\Omega= {\mathbb  R}^n $ for
${\Delta}$, ${\Delta}_{D}$,  ${\Delta}_{N}$,
${\Delta}_{DN}$; $\Omega= {\mathbb  R}^n_+ $ for ${\Delta}_{N_+}$,
${\Delta}_{D_+}$
and $\Omega= {\mathbb  R}^n_- $ for  ${\Delta}_{N_-}$,  ${\Delta}_{D_-}$.

      \medskip

($\gamma$) If $L$ is  one of the operators  ${\Delta}$, ${\Delta}_{N_+},$
      ${\Delta}_{N_-}$ and   ${\Delta}_{N}$, then $L$ conserves probability,
that is
$$
\exp(-tL)\one=\one.
$$
This conservative property   does not hold for  ${\Delta}_D$, ${\Delta}_{D_+},$
      ${\Delta}_{D_-}$  and  ${\Delta}_{DN}$.

\subsection{{\rm BMO} spaces associated with operators}\label{sec21}

Suppose that $\Omega \subset {\mathbb  R}^n$ is an open subset of
${\mathbb  R}^n$.
Suppose that $L$ is a linear operator on $L^2(\Omega)$
which
generates an analytic semigroup $e^{-tL}$ with a  kernel $p_t(x,y)$  satisfying
Gaussian
upper bound (\ref{e2.1}).

We define
$$
{\mathcal M}(\Omega) =\bigg\{f\in L^1_{\rm loc}({\Omega}): \exists d>0, \ \
\int_{\Omega}\frac{ |f(x)|^2}{   1+|x |^{n+d}  } dx<\infty
\bigg\}.
$$
Note that in virtue of the Gaussian bounds (\ref{e2.1}) we can extend the
action of the semigroup
operators  $\exp(-tL)$ to the space ${\mathcal M}(\Omega)$, that is
we can define
$\exp(-tL)f$ for all  $f\in{\mathcal M}(\Omega)$.
By $B(x,r)$ we denote the ball in $\Omega$ with respect to the Euclidean
distance restricted to $\Omega$ that is
$B(x,r)=  \{y\in \Omega \colon \, |x-y| < r \}   $.

The following {\BLO} space associated with
an operator $L$  was introduced in \cite{DY1}.

\begin{definition}
\label{defi23}
We say that $f\in {\mathcal M}(\Omega)$ is of bounded  mean oscillation
associated with an operator $L$  (abbreviated as
{\BLO}) if
\begin{eqnarray}
\label{e2.3}
      \|f\|_{{\rm BMO}_L(\Omega)}=\sup_{B(y,r)}\frac{1}{|B(y,r)|}\int_{B(y,r)}|
f(x)-\exp(-{r^2}L)f(x)|dx <\infty,
\end{eqnarray}
where the supremum  is taken over all balls $B(y,r)$  in $\Omega$.

The smallest bound  for which {\rm (\ref{e2.3})} is satisfied is then taken to
be the norm of $f$ in this space, and is denoted by
$\|f\|_{{\rm BMO}_L(\Omega)}$.
\end{definition}

\noindent
{\bf Remarks.}
{\rm (i)}\
Note that ({\BLO},
$\|\cdot\|_{{\rm BMO}_L(\Omega)}$) is a semi-normed
vector space, with the  semi-norm vanishing on the  kernel space
${\mathcal K}_L$  defined by
\begin{eqnarray*}
{\mathcal K}_L=\Big\{f\in{\mathcal M}(\Omega)\colon \,
\exp(-tL)f=f, \ \      \quad \forall t>0 \Big\}.
\end{eqnarray*}
The class of functions of ${\rm BMO}_L(\Omega)$ (modulo
${\mathcal K}_L$)
is a Banach  space.
We refer the reader to Section 6 of \cite{DY2} for a  discussion   on   the
dimension   of the   space  ${\mathcal K}_L$ of
${\rm BMO}_L({\mathbb R}^n)$ when $L$ is
      a second order divergence form elliptic operator
or a Schr\"odinger operator.
In the sequel By ${\rm BMO}_L(\Omega)$ we always denote the space
${\rm BMO}_L(\Omega)$ (modulo
${\mathcal K}_L$) and we skip (modulo
${\mathcal K}_L$) to simplify notation.

\medskip

{\rm (ii)}\ Similarly to the classical BMO space, it is easy to check that
$L^{\infty}(\Omega) \subset  {\rm BMO}_L(\Omega)$  with
$\|f\|_{{\rm BMO}_L(\Omega)}\leq 2\|f\|_{L^\infty}.$

\medskip

{\rm (iii)}\
The classical
BMO  space (modulo all constant functions)  and the
$ {\rm BMO}_{{\Delta}}({\mathbb R}^n)$ space (modulo all
       harmonic  functions),
       coincide, and their norms are equivalent.
See Theorem~2.15 of \cite{DY1}.

\medskip

{\rm (iv)}\
Note that the Euclidean distance in Definition~\ref{defi23} can be replaced
by any equivalent distance. That is if there exists $c>0$ such that
$c^{-1}|x-y| \le d(x,y) \le c|x-y|$ then one can take in (\ref{e2.3}) the
supremum over
all balls $B^d(x,r)$ with respect to the metric $d$.  In particular if
$\Omega={\mathbb R}^n$,
$\Omega={\mathbb R}^n_+$ or
$\Omega={\mathbb R}^n_-$, one can take the supremum over all cubes
$Q$ such that
$Q \subset \Omega$ in (\ref{e2.3}), i.e.,
   we  can define equivalent norm in {\BLO}
by the formula
\begin{eqnarray}
\label{ecube}
      \|f\|_{{\br}_L({\mathbb R}^n)}=\sup_{Q}\frac{1}{|Q|}\int_{Q}
|f(x)-\exp(-{l_Q^2}L)f(x)|dx <\infty,
\end{eqnarray}
where $l_Q$
is the side length of $Q$ and the supremum is taken over all
cubes $Q\subset \Omega$.

\bigskip

The following proposition is essentially equivalent to
Proposition~3.1 of \cite{Ma}.

\begin{prop}\label{prop24} Assume that   for every $\ \! t>0$,
$e^{-tL}(\one)=\one$ almost everywhere, that is,
$
\int_{{\mathbb R}^n} p_t(x,y)dy=1\ {\rm for}\  {\rm almost\ all}\
x\in{\mathbb R}^n.
$
Then, we have ${\rm BMO}({\mathbb R}^n)\subset{\rm BMO}_L({\mathbb R}^n)$,
and there exists a
positive constant
$c>0$ such that
\begin{equation}
\|f\|_{{\rm BMO}_L({\mathbb R}^n)}\leq c\|f\|_{{\rm BMO}({\mathbb R}^n)}.
\label{e2.5}
\end{equation}
However, the converse inequality does not hold in general.
\end{prop}

We remark that condition $e^{-tL}(\one=\one,$  is necessary for
(\ref{e2.5}). Indeed,  (\ref{e2.5})
implies $\|\one\|_{{\rm BMO}_L({\mathbb R}^n)}=0$. Hence
$e^{-tL}(\one)=\one$ almost everywhere for all $t>0$,

\bigskip

\section{BMO spaces on the  half spaces and
BMO spaces associated with the Dirichlet and Neumann Laplacian.}\label{half}
\setcounter{equation}{0}


In this section we  describe the equivalence
   between the {\BR} spaces on the half space and
{\BR}  spaces corresponding
to the Neumann and Dirichlet Laplacian.

\begin{prop}\label{prop22}
    {\rm (i)} The spaces ${\rm BMO}_r({\mathbb
R}^n_+)$  and
${\rm BMO}_e({\mathbb R}^n_+)$   coincide, and their
norms are equivalent.

\smallskip

{\rm (ii)} The spaces   ${\rm BMO}_z({\mathbb R}^n_+)$  and
     ${\rm BMO}_o({\mathbb R}^n_+)$  coincide, and their
norms are equivalent.
\end{prop}

\begin{proof}
Following  \cite{CKS}, for any function $f \in L^1({\mathbb R}^n_+)$
we set
\begin{eqnarray}
\|f\|_{H^1_e({\mathbb R}^n_+)} =\|f_e\|_{H^1({\mathbb R}^n)}
      \quad \mbox{and} \quad
      \|f\|_{H^1_o({\mathbb R}^n_+)} =\|f_o\|_{H^1({\mathbb R}^n)}
\label{e2.9}
\end{eqnarray}
    and by ${H^1_e({\mathbb R}^n_+)}$  and ${H^1_o({\mathbb R}^n_+)}$
     we denote the corresponding   Banach spaces.
It follows from  Corollaries 1.6, 1.8 of \cite{CKS} and Proposition
32 of \cite{ART} that
     the dual space of ${H^1_e({\mathbb R}^n_+)}$
is  the space ${{\rm BMO}_r({\mathbb R}^n_+)}$  and
the dual space of ${H^1_o({\mathbb R}^n_+)}$
is  the space ${{\rm BMO}_z({\mathbb R}^n_+)}$ . See also \cite{AR}.

 The inclusion
$ {\rm BMO}_e({\mathbb R}^n_+)\subseteq {\rm BMO}_r({\mathbb R}^n_+)
$  is obvious. Hence to  prove (i) it is enough to show that
  $ {\rm BMO}_r({\mathbb R}^n_+)\subseteq
{\rm BMO}_e({\mathbb
R}^n_+).$ Let
$f\in {\rm BMO}_r({\mathbb R}^n_+)$.  To see that $f\in {\rm
BMO_e({\mathbb R}^n_+)},$
by the definition it reduces to proving
      $f_e\in {\rm BMO({\mathbb R}^n)} $ where $f_e$ is the even
extension of $f$. For any $g(x)\in H^1({\mathbb R}^n)$, we denote
by ${\widetilde g}(x)=g({\widetilde x})$ where ${\widetilde x}=(x',
-x_n)$. Since  $({H^1_e({\mathbb R}^n_+)})'={{\rm
BMO}_r({\mathbb
R}^n_+)}$,  we   have
\begin{eqnarray*}
|\int_{{\mathbb R}^n} f_e(x)g(x)dx|&= &|\int_{{\mathbb R}^n_-} f_e(x)g(x)dx
+\int_{{\mathbb R}^n_+} f_e(x)g(x)dx|\\
&=&|\int_{{\mathbb R}^n_+} f(x) \Big( g({\widetilde  x})+g(x)\Big)dx|\\
&\leq& c\|f\|_{{\rm BMO}_r({\mathbb R}^n_+) } \|({\widetilde
g}+g)\|_{H^1_e({\mathbb
R}^n_+)}\\
&\leq &c\|f\|_{{\rm BMO}_r({\mathbb R}^n_+) } \|g\|_{H^1({\mathbb
R}^n)}.
\end{eqnarray*}
This shows that
$ {\rm BMO}_r({\mathbb R}^n_+)\subset {\rm
BMO_e({\mathbb R}^n_+)}$, and proves {\rm (i)}.

We now prove {\rm (ii)}. The inclusion
     $ {\rm BMO}_z({\mathbb R}^n_+)\subseteq {\rm BMO}_o({\mathbb R}^n_+)
$ is  obvious.
Let
$f\in {\rm BMO}_o({\mathbb R}^n_+)$ and thus  $f_o\in {\rm
BMO({\mathbb R}^n)}.$
To see that $f\in {\rm BMO}_z({\mathbb R}^n_+)$,
it reduces to proving $f\in ({H^1_o({\mathbb R}^n_+)})'$
since
${\rm BMO}_z({\mathbb R}^n_+)=({H^1_o({\mathbb R}^n_+)})'$.
If
$g\in H^1_o({\mathbb R}^n_+)$, then  $g_o\in H^1({\mathbb R}^n).$
Hence
\begin{eqnarray*}
    |\int_{{\mathbb R}^n_+} f(x)g(x)dx|&=&
    \frac{1}{2}|\int_{{\mathbb R}^n} f_o(x)g_o(x)dx|\\
&\leq&  c\|f_o\|_{{\rm BMO}({\mathbb R}^n) }\|g_o\|_{H^1({\mathbb
R}^n)}\\
&\leq&  c\|f\|_{{\rm BMO}({\mathbb R}^n_+) }\|g\|_{H^1({\mathbb
R}^n_+)}.
\end{eqnarray*}
This shows that
$ {\rm BMO}_o({\mathbb R}^n_+)\subset {\rm
BMO_z({\mathbb R}^n_+)}$, and proves {\rm (ii)}.
\end{proof}
We use Proposition~\ref{prop22} to obtain the following result.
\begin{prop}
\label{prop26}   {\rm (i)} The spaces
${\rm BMO}_{{\Delta}_{D_+}}({\mathbb R}^n_+)$,
      ${\rm BMO}_z({\mathbb R}^n_+)$  and  ${\rm BMO}_o({\mathbb R}^n_+)$
      coincide, and their norms are equivalent.

\smallskip

{\rm (ii)} The spaces
${\rm BMO}_{{\Delta}_{N_+}}({\mathbb R}^n_+)$,
      ${\rm BMO}_r({\mathbb R}^n_+)$  and  ${\rm BMO}_e({\mathbb R}^n_+)$
      coincide, and their norms are equivalent.
\end{prop}
\begin{proof} We first prove {\rm (i)}.
Let $f\in
{\rm BMO}_z({\mathbb R}^n_+)$. By  Proposition~\ref{prop22} we have that
$f\in {\rm BMO}_o({\mathbb R}^n_+)$ and then $f_o\in {\rm BMO}({\mathbb R}^n)$.
To
prove $f\in
      {\rm BMO}_{{\Delta}_{D_+}}({\mathbb R}^n_+)$, it suffices to show that
for any cube
      $Q\subseteq {\mathbb R}^n_+$,
\begin{eqnarray}
      \int_Q|f(x)-e^{-{l^2_Q}{\Delta}_{D_+}}f(x)|dx\leq c|Q|\|f\|_{{\rm
BMO}_z({\mathbb R}^n_+) }.
\label{e2.15}
\end{eqnarray}
By (\ref{e2.5}) and Propositions~\ref{prop22}
\begin{eqnarray*}
\frac{1}{|Q|}
\int_Q|f(x)-e^{-{l^2_Q}{\Delta}_{D_+}}f(x)|dx&=&\frac{1}{|Q|}
\int_{Q}|f(x)-e^{-l^2_Q{\Delta}} f_o(x)|dx \\
&\leq&  c\|f_o\|_{{\rm BMO}({\mathbb R}^n)}
\leq  c\|f\|_{{\rm BMO}_o({\mathbb R}^n_+)}\\
&\leq&  c\|f\|_{{\rm BMO}_z({\mathbb R}^n_+) }.
\end{eqnarray*}
This proves (\ref{e2.15}).

Next assume that   $f\in {\rm BMO}_{{\Delta}_{D_+}}({\mathbb R}^n_+)$.
By Proposition~\ref{prop22},
$f\in  {\rm BMO}_z({\mathbb R}^n_+)$ or equivalently
      $f_o\in {\rm BMO}({\mathbb R}^n)$.
Note that by (\ref{e2.12})
it is enough  to  prove that
for any cube
      $Q\subseteq {\mathbb R}^n$,
\begin{eqnarray}
\int_Q|f_o(x)-e^{-l^2_Q{\Delta}} f_o(x)|dy\leq c
|Q|\|f\|_{ {\rm BMO}_{{\Delta}_{D_+}}({\mathbb R}^n_+)}.
\label{e2.16}
\end{eqnarray}
We now   verify (\ref{e2.16}). Let us examine the  cubes
     $Q.$

\smallskip
Case 1:\   If $Q\subseteq {\mathbb R}^n_-$, then
     for any $x\in Q,$
$$
-\exp(-{l^2_Q}{\Delta}_{D_+})f({\widetilde x})=
\exp(-{l^2_Q}{\Delta}) f_o(x)
$$
and ${\widetilde x}\in \widehat{Q}\subseteq{\mathbb R}^n_+$
(here   $\widetilde{Q}$ is   a cube
defined   in (\ref{e2.7})). Note also that $|\widetilde{Q}|=|Q|$. Hence
\begin{eqnarray*}
\int_Q|f_o(x)-e^{-l^2_Q{\Delta}} f_o(x)|dx
= \int_{\widetilde{Q}}|f_o({\widetilde x})-
e^{-{l^2_{\widetilde{Q}}}{\Delta}_{D_+}}f({\widetilde x})|dx
\leq c|Q|\|f\|_{ {\rm BMO}_{{\Delta}_{D_+}}({\mathbb R}^n_+)}.
\end{eqnarray*}

Case 2:\ If $Q\cap {\mathbb R}^n_-\not=\emptyset$ and
$Q\cap {\mathbb R}^n_+\not=\emptyset$, then let  ${{\widehat Q}_-}$
and ${{\widehat Q}_+}$ be the two cubes as in (\ref{e2.8}).
By (\ref{e2.12})   and Proposition~\ref{prop22},
\begin{eqnarray*}
      \int_Q|f_o(x)-e^{-l^2_Q{\Delta}} f_o(x)|dx&=&
      \int_{Q_-\cup Q_+}|f_o(x)-e^{-l^2_Q{\Delta}} f_o(x)|dx\\
&\leq&
2\int_{{\widehat Q}_+}|f(x)-   e^{-{l^2_{{Q}}}{\Delta}_{D_+}}f(x)|dx \\
&\leq&  2|Q|\|f\|_{{\rm BMO}_{{\Delta}_{D_+}}({\mathbb R}^n_+)}.
\end{eqnarray*}

Case 3:\ If $Q\subseteq {\mathbb R}^n_+$,
then
$e^{-l^2_Q{\Delta}} f_o(x)=  e^{-l^2_Q{\Delta}_{D_+}} f(x)$  for any $x\in Q$.
Hence
\begin{eqnarray*}
\int_Q|f_o(x)-e^{-l^2_Q{\Delta}} f_o(x)|dx\leq |Q|
\|f\|_{ {\rm BMO}_{{\Delta}_{D_+}}({\mathbb R}^n_+)}.
\end{eqnarray*}
The estimate (\ref{e2.16}) follows readily.
This shows that
$f_o\in {\rm BMO}({\mathbb R}^n)$ so    $f\in
{\rm BMO}_z({\mathbb R}^n_+)$.

The proof of {\rm (ii)} is similar to the proof of {\rm (i)}  so we skip it.
\end{proof}

In a similar way  as for the upper-half space,   we can  define the
space ${\rm
BMO}_{{\Delta}_{D_-}}({\mathbb R}^n_-)$ and
${\rm BMO}_{{\Delta}_{N_-}}({\mathbb R}^n_-)$
associated with the  Dirichlet and  Neumann
Laplacian ${{\Delta}_{D_-}}$, ${{\Delta}_{N_-}}$
on the lower-half space ${\mathbb R}^n_-$.
The same argument as in Proposition~\ref{prop26} gives the following
proposition. We leave  the proof to the reader.

\begin{prop}
\label{prop28} {\rm (i)}\ The spaces ${\rm BMO}_{{\Delta}_{D_-}}
({\mathbb R}^n_-)$,
      ${\rm BMO}_z({\mathbb R}^n_-)$  and ${\rm BMO}_o({\mathbb R}^n_-)$
      coincide, and their norms are equivalent.

{\rm (ii)}\ The  spaces ${\rm BMO}_{{\Delta}_{N_-}}({\mathbb R}^n_-)$,
${\rm BMO}_r({\mathbb R}^n_-)$ and
${\rm BMO}_e({\mathbb R}^n_-)$   coincide, and their norms are equivalent.
\end{prop}

\bigskip

\section{Comparison  between   the classical BMO  and   the  new  BMO spaces
      associated  with operators}\label{comp}
\setcounter{equation}{0}

In the introduction we mention that
all cases of relation between  the classical BMO  and   the  new  BMO
spaces are possible. The following theorem provides simple example to
prove this
statement.

\begin{theorem}\label{theo41}
The following inclusions describe the relation between
${\rm BMO}_{{\Delta}_{D}}({\mathbb R}^n) $,
${\rm BMO}({\mathbb R}^n)$ and
${\rm BMO}_{{{\Delta}_{N}}}({\mathbb R}^n)$
\begin{eqnarray}
{\rm BMO}_{{\Delta}_{D}}({\mathbb R}^n) \subsetneqq
{\rm BMO}({\mathbb R}^n)
\subsetneqq {\rm BMO}_{{{\Delta}_{N}}}({\mathbb R}^n).
\label{e4.2}
\end{eqnarray}
That is, the classical ${\rm BMO}$ space is a proper subspace of ${\rm
BMO}_{{\Delta}_{N}}({\mathbb R}^n)$, and ${\rm
BMO}_{{\Delta}_{D}}({\mathbb R}^n)$ is a proper subspace of
    ${\rm BMO}$.

Moreover, we have
\begin{eqnarray}
{\rm BMO}({\mathbb R}^n)\not\subseteq {\rm
BMO}_{{\Delta}_{DN}}({\mathbb R}^n)
\quad \mbox{ and} \quad
{\rm BMO}_{{\Delta}_{DN}}({\mathbb R}^n)\not\subseteq {\rm
BMO}({\mathbb R}^n).
\label{e4.3}
\end{eqnarray}
\end{theorem}
    The proof of Theorem~\ref{theo41} is based on the following proposition.
\begin{prop}\label{prop42}
The {\rm BMO} spaces corresponding to the operators
${{\Delta}_{N}}$
${{\Delta}_{D}}$
${{\Delta}_{ND}}$ can be described in the following way
\begin{eqnarray*}
{\rm BMO}_{{\Delta}_{N}}({\mathbb R}^n)= \Big\{f\in {\mathcal
M}({\mathbb R}^n): \
      f_+ \in {\rm BMO}_r({\mathbb R}^n_+) \ {\rm and} \
f_-  \in {\rm BMO}_r({\mathbb R}^n_-)\Big\};
\\
{\rm BMO}_{{\Delta}_{D}}({\mathbb R}^n)= \Big\{f\in {\mathcal
M}({\mathbb R}^n): \
      f_+ \in {\rm BMO}_z({\mathbb R}^n_+) \ {\rm and} \
f_-  \in {\rm BMO}_z({\mathbb R}^n_-)\Big\};
\\
{\rm BMO}_{{\Delta}_{DN}}({\mathbb R}^n)= \Big\{f\in {\mathcal
M}({\mathbb R}^n): \
      f_+ \in {\rm BMO}_r({\mathbb R}^n_+) \ {\rm and} \
f_-  \in {\rm BMO}_z({\mathbb R}^n_-)\Big\}.
\end{eqnarray*}
\end{prop}
\begin{proof}
In the following proof $L$ is one of the operators  ${\Delta}_{N}$,
${\Delta}_{D}$
or ${\Delta}_{DN}$. If $L={\Delta}_{N}$, then we denote by
$L_+={\Delta}_{N_+}$ and $L_-={\Delta}_{N_-}$. Similarly if
$L={\Delta}_{D}$ then
$L_+={\Delta}_{D_+}$ and $L_-={\Delta}_{D_-}$. Finally for
$L={\Delta}_{DN}$ we let
$L_+={\Delta}_{N_+}$ and $L_-={\Delta}_{D_-}$.

By (\ref{e2.19n}), (\ref{e2.19d}) and (\ref{e2.19dn})
\begin{equation}\label{e4.19}
(\exp(-tL) f)_+=\exp(-t{L}_{+})f_+
\quad \mbox{and} \quad
(\exp(-tL) f)_-=\exp(-t{L}_{-})f_-
\end{equation}
for any of the three considered operators.
Hence for any cube $Q\subset {\mathbb
R}^n$ we have
\begin{eqnarray}
\int_{Q}|f-e^{-{l^2_Q}{L}} f(x)|dx&=&\int_{Q\cap {
{{\mathbb R}^n_-}}}|f_- - e^{-{l^2_Q}{{L_-}}} f_-(x)|dx\nonumber \\ &+&
\int_{Q\cap { {{\mathbb R}^n_+}}}|f_+-   e^{-{l^2_Q}{{L_+}}}f_+(x)|dx.
\label{e2.21}
\end{eqnarray}
In virtue of Propositions ~\ref{prop26} and ~\ref{prop28} it is
enough to show that
$$
{\rm BMO}_{L}({\mathbb R}^n)= \Big\{f\in {\mathcal M}({\mathbb R}^n): \
      f_+ \in {\rm BMO}_{L_+}({\mathbb R}^n_+) \ {\rm and} \
f_-  \in {\rm BMO}_{L_-}({\mathbb R}^n_-)\Big\}.
$$
      Assume now that
$f\in {\mathcal M}({\mathbb R}^n)$  such that
$f_-\in {\rm BMO}_{L_-}({\mathbb R}^n_-)$ and
$f_+\in {\rm BMO}_{L_+}({\mathbb R}^n_+)
$.
In order to prove
      $f \in {\rm BMO}_{L}({\mathbb R}^n)$, it suffices to prove
       that
for any cube
      $Q\subseteq {\mathbb R}^n$,
\begin{eqnarray*}
\int_Q|f(x)-e^{-{l^2_Q}L} f (x)|dy\leq c
|Q|\Big(\|f_-\|_{ {\rm BMO}_{L_-}({\mathbb R}^n_-)} +
      \|f_+\|_{ {\rm BMO}_{L_+}({\mathbb R}^n_+)}\Big).
\end{eqnarray*}
As in the proof of Proposition~\ref{prop26}, we consider the
following three cases of $Q$.

\smallskip

Case 1: \  If $Q\subseteq {\mathbb R}^n_-$, then by (\ref{e2.21})
\begin{eqnarray*}
      \int_Q|f(x)-e^{-{l^2_Q} {L}} f(x)|dx&=&
\int_Q|f_-(x)-e^{-l^2_Q L_-}f_-(x)|dx \\
      &\leq& c|Q|\|f_-\|_{{\rm BMO}_{L_-}({\mathbb R}^n_-)}.
\end{eqnarray*}

Case 2: \ If  $Q\cap {\mathbb R}^n_-\not=\emptyset$ and
$Q\cap {\mathbb R}^n_+\not=\emptyset$, then let  ${\widetilde {Q}_-}$
and ${\widetilde {Q}_+}$ be the cubes as in (\ref{e2.8}).
By (\ref{e2.21})
\begin{eqnarray*}
\int_Q|f(x)-&&\hspace{-1cm}e^{-{l^2_Q}{L}}f(x)|dx=
\int_{Q_-\cup Q_+}|f(x)-e^{-l^2_Q L}f(x) |dx \\
&\leq&  \int_{\widetilde {Q}_-}|f_-(x)-  e^{-{l^2_Q} L_- }  f_-(x)  |dx
    + \int_{\widetilde {Q}_+}|f_+(x)- e^{-{l^2_{{Q}}} {{L_+}}}f_+(x)  |dx\\
&\leq&c|Q|\Big(\|f_-\|_{{\rm BMO}_{L_-}({\mathbb R}^n_-) }
+ \|f_+\|_{{\rm BMO}_{L_+}({\mathbb R}^n_+) }\Big).
\end{eqnarray*}

Case 3:\  If $Q\subseteq {\mathbb R}^n_-$, then by
(\ref{e2.21})
\begin{eqnarray*}
      \int_Q|f(x)-  e^{-{l^2_Q}L} f(x)|dx
&=& \int_Q|f_+(x)-  e^{-{l^2_{{Q}}} { {L_+}}} f_+(x)|dx\\
       &\leq& c|Q|\|f_+\|_{{\rm BMO}_{\Delta_{N_+}}({\mathbb R}^n_+)}.
\end{eqnarray*}
Hence $f\in {\rm BMO}_{L}({\mathbb R}^n).$

We now assume that $f\in {\rm BMO}_{L}({\mathbb R}^n)$.
     By (\ref{e2.21}), we have  that $f_-\in {\rm
BMO}_{{L_-}}({\mathbb R}^n_-)$
      and   $f_+\in {\rm BMO}_{\Delta_{N_+}}({\mathbb R}^n_+)$.
Now Proposition~\ref{prop42} is a straightforward
consequence of Propositions ~\ref{prop26} and ~\ref{prop28}.
\end{proof}

The logarithmic function is  a simple example that typifies some of
the essential properties of the classical space BMO.
For example if we define  function $\log \colon \, {\mathbb R}^n \to
{\mathbb R}$
by the formula $\log^e(x)=\log|x_n|$ for all $x\in {\mathbb R}^n$ and
${\rm Log}(x)=H(x_n)\log|x_n|$, where $H$ is the Heaviside function then
\begin{eqnarray}
\begin{array}{l}
      {\log^e} \in {\rm BMO}({\mathbb R}^n)
\\
     {\rm Log}
\not\in {\rm BMO}({\mathbb
R}^n).
\end{array}
\label{log}
\end{eqnarray}
See, for examples,
Chapter IV of \cite{St2} and page 217   of \cite{To} .
    We will use the property (\ref{log}) in the proof of
Theorem~\ref{theo41}

\begin{proof}[Proof of Theorem~\ref{theo41}]
It is a straightforward consequence of    Definition~\ref{defi21}  that
if $f_+ \in  {\rm BMO}_z({\mathbb R}^n_+)$
and $f_- \in  {\rm BMO}_z({\mathbb R}^n_-)$
then $f \in {\rm BMO}$. It also  follows from Definition~\ref{defi21}
that if $f \in {\rm BMO}$ then $f_+ \in  {\rm BMO}_r({\mathbb R}^n_+)$
and  $f_- \in  {\rm BMO}_r({\mathbb R}^n_-)$.
Hence it follows from Theorem~\ref{theo41} and
Propositions~\ref{prop26}~and~\ref{prop28} that
$$
{\rm BMO}_{{\Delta}_{D}}({\mathbb R}^n) \subset
{\rm BMO}({\mathbb R}^n)
\subset {\rm BMO}_{{{\Delta}_{N}}}({\mathbb R}^n).
$$
To prove that the above inclusions are proper  we note
that by (\ref{log}) and Definition~\ref{defi21}
$$
\log_+\notin {\rm BMO}_z({\mathbb R}^n_+)
\quad \mbox{and} \quad \log_+\in {\rm BMO}_r({\mathbb R}^n_+),
$$
    where $\log_+$ is
the restriction of $\log^e$ to ${\mathbb R}^n_+$. Next if $\log_-$ is
the restriction of $\log^e$ to ${\mathbb R}^n_-$ then
$$
\log_-\notin {\rm BMO}_z({\mathbb R}^n_-) \quad
\mbox{and} \quad \log_-\in {\rm BMO}_r({\mathbb R}^n_-).
$$
Hence
$$
\log^e \in  {\rm BMO} \quad \mbox{and} \quad
{\log^e}   \notin {\rm BMO}_{{\Delta}_{D}}({\mathbb R}^n).
$$
Similarly
$$
{\rm Log} \notin  {\rm BMO}\quad \mbox{and} \quad
{\rm Log}   \in {\rm BMO}_{{\Delta}_{N}}({\mathbb R}^n)
$$
This ends the proof of (\ref{e4.2}).
Finally to prove (\ref{e4.3}) we note that ${\rm Log}
\in {\rm BMO}_{{\Delta}_{DN}}({\mathbb R}^n)$
and $\log   \notin {\rm BMO}_{{\Delta}_{DN}}({\mathbb R}^n)$.
\end{proof}

    \noindent
{\bf Remark.}\
Suppose that  $L$ is  a
linear operator  on
$L^2({\mathbb R}^n)$ which generates
      an analytic semigroup $e^{-tL}$ with kernels
$p_t(x,y)$  satisfying   upper bound (\ref{e2.1}). Under the
additional  condition that  the kernel
$p_t(x,y)$ of
$e^{-tL}$  has sufficient regularities on space  variables
$x, y$ and
$e^{-tL}(\one)=e^{-tL^{\ast}}(\one)=\one$,  it can be proved that   classical space
BMO   and  the space
$ {\rm BMO}_L({\mathbb R}^n)$  spaces  coincide, and their norms
are equivalent. See Section 6 of \cite{DY2}.

\medskip

Next we discuss
    the duality of the Hardy   and BMO spaces associated with operators.
Suppose that  $L$ is  a
linear operator  on
$L^2({\mathbb R}^n)$ which generates
      an analytic semigroup $e^{-tL}$ with kernels
$p_t(x,y)$  satisfying   Gaussian  upper bound (\ref{e2.1}).
For any $(x,t)\in {\mathbb R}^n\times (0, \infty)$, we define
$$
Q_{t}f(x)= -t\frac{d}{ dt}e^{-tL}f(x)=tLe^{-tL}f(x)
$$
for any $f\in{\mathcal M}$. Following \cite{ADM}, given a function
$f\in L^1({\mathbb
R}^n)$, the area integral function ${\mathcal S }_{L}(f)$ associated
with an operator
$L$ is defined by
$$
{\mathcal S }_{L}f(x)=\Big( \int_0^{\infty}\!\!\!\!\int_{|y-x|<t}\big|
Q_{t^2}f(y)\big|^2\
\frac{dy \ \! dt}{t^{n+1}} \Big)^{1/2}, \ \ \ \ \  x\in {\mathbb R}^n.
$$
The following definition was introduced in \cite{ADM}.    We say that
$f\in L^1({\mathbb R}^n)$   belongs   to a  Hardy space associated
with $L$ (abbreviated as
$H^1_{L}({\mathbb R}^n))$ if
$
{\mathcal S }_{L}f \in L^1.
$
If it is the case, we define its norm by
$$
\|f\|_{H^{ 1}_{L}({\mathbb R}^n)}=\|{\mathcal S }_{L}f
\|_{L^1}.
$$
Note that if $L={\Delta}$ is the Laplacian on ${\mathbb R}^n$, then
the classical
Hardy space $H^1$ and
$ H^1_{{\Delta}}$
       coincide, and their norms are equivalent. See \cite{ADM}.

Under  the   assumptions  that   $L$ satisfies  Gaussian upper
bound (\ref{e2.1})
and has a bounded
$H_{\infty}$-calculus in $L^2({\mathbb R}^n)$,  it was proved in
\cite{DY2} that the dual space of the $H^1_L({\mathbb R}^n)$  space is the ${\rm
BMO}_{L^{\ast}}
({\mathbb R}^n)$ space in which $L^{\ast}$ is the adjoint operator of $L$.

\smallskip

Note that the operators ${\Delta}_D$, ${\Delta}_N$ and
${\Delta}_{DN}$ are self-adjoint operators, hence each of them has a bounded
$H_{\infty}$-calculus in $L^2({\mathbb R}^n)$. See \cite{Mc}.
We thus have the following
corollary.

\begin{cor}\label{cor37}
{\rm (i)}\ The dual space of
$H^1_{{\Delta}}({\mathbb R}^n)$  is the space
${\rm BMO}_{{\Delta}}({\mathbb R}^n).$

\medskip

{\rm (ii)}\ The dual spaces of
$H^1_{{\Delta}_D}({\mathbb R}^n)$, $H^1_{{\Delta}_N}({\mathbb R}^n)$
or $H^1_{{\Delta}_{DN}}({\mathbb R}^n)$
are the spaces ${\rm BMO}_{{\Delta}_D}({\mathbb
R}^n) $, ${\rm BMO}_{{\Delta}_N}({\mathbb R}^n)$ or
${\rm BMO}_{{\Delta}_{DN}}({\mathbb R}^n)$, respectively.

\medskip

{\rm (iii)}\ For the Neumann Laplacian ${\Delta}_N$ on
${\mathbb R}^n$, we have that
$H^1_{{\Delta}_N}({\mathbb R}^n) \subsetneqq
      H^1({\mathbb R}^n)$ and $H^1_{{\Delta}_N}({\mathbb R}^n)\not=\emptyset$.
That is, $  H^1_{{\Delta}_N}({\mathbb R}^n)$ is a proper subspace of
the classical Hardy space
$H^1({\mathbb R}^n)$.
\end{cor}

\noindent
{\bf Remark.} In \cite{W}, it was asked if a proper subspace of the classical
Hardy space exists in which the subspace is characterized by maximal
functions. This
question was answered positively in  \cite{UW}. Our result
     (iii) of Corollary~\ref{cor37} gives a proper subspace
of the classical
Hardy space where  the subspace is characterized by area integral  functions.

\bigskip

\section{ Fractional powers  $L^{-\alpha/2}$ and the space
${\rm BMO}_L({\mathbb R}^n)$}\label{frac}
\setcounter{equation}{0}

\subsection{ Boundedness of fractional powers  $L^{-\alpha/2}$.  }
For any $0<\alpha<n, $       the
fractional powers
$L^{-\alpha/2}$ of
$L$   is defined by
\begin{eqnarray}
L^{-\alpha/2}f(x) =\frac{1}{\Gamma(\alpha/2)}\int_0^{\infty}
    t^{{\alpha/2} -1}  e^{-tL}f(x){dt}.
\label{e3.1}
\end{eqnarray}
We assume that the  semigroup $e^{-tL}$ has a  kernel
$p_t(x,y)$ which  satisfies the
upper bound (\ref{e2.1}) so
$
|L^{-\alpha/2}f(x)|\leq c{\mathcal I}_{\alpha}(|f|)(x)
$
      for all
$
x\in{\mathbb R}^n$, where
$$
{\mathcal I}_{\alpha}f(x)=\int_{{\mathbb R}^n} \frac{f(y)}{
|x-y|^{n-\alpha} }dy, \ \ \ \ \ 0<\alpha<n,
$$
      is the classical
fractional powers  of the Laplacian ${\Delta}$ on ${\mathbb R}^n$.

Let us recall that the semigroup $\{\exp(-tL)\colon \, t>0\}$
acting on $L^p({\mathbb R}^n)$
is equicontinuous on  $L^p ({\mathbb R}^n)  $ if
$\sup_{t>0} \|e^{-tL}\|_{L^p \to L^p} < \infty$. Note that all the
semigroups which
we
consider here are equicontinuous on all $L^p({\mathbb R}^n)$
for $1\leq p \leq \infty$.
In the sequel we  need the  following Hardy-Littlewood-Sobolev
theorem.  See Theorem~II.2.7 page 12 \cite{VSC}.

\begin{prop}\label{prop41}
Suppose that $e^{-tL}$ is a semigroup
which is equicontinuous on $L^1({\mathbb R}^n)$ and
$L^{\infty}({\mathbb R}^n)$.
Also suppose that
$$
p_t(x,x) \le t^{-n/2}.
$$
Then for $0<\alpha<n$,

\medskip

{\rm (i)} \  for $1< p<\frac{n}{\alpha}$ and $ \frac{1}{ q}=\frac{1}{p}
-\frac{\alpha}{n}$, we
   have
$$\|L^{-\alpha/2}f\|_{L^q}\leq c_{p,q}\|f\|_{L^p};
$$

{\rm (ii)}\ $L^{-\alpha/2}$ is of weak-type $(1,q)$, that is, for any
$\lambda>0$, we have
$$
\Big|\{x: |L^{-\alpha/2}f(x)|>\lambda\}\Big|\leq c
\Big(\frac{\|f\|_{L^1}}{ \lambda}\Big)^q,
$$
where $ q=(1-\frac{\alpha}{n})^{-1}.$
\end{prop}

\bigskip

Let us consider  the limiting case  $q=\infty$ in Proposition~\ref{prop41}.
It is well-known that for every $f\in L^{n/\alpha}({\mathbb R}^n)$,
either ${\mathcal
I}_{\alpha}f
\equiv \infty$ or ${\mathcal I}_{\alpha}f
\in {\rm BMO}({\mathbb R}^n)$ with
\begin{equation}\label{toto}
\|{\mathcal I}_{\alpha}f \|_{{\rm BMO}({\mathbb R}^n)}
\leq c\|f\|_{L^{n/\alpha}},
\end{equation}
see
page 221 of \cite{To}. An example of  ${\mathcal I}_{\alpha}f\equiv \infty $
is given by $f(x)=
|x|^{-\alpha} {\rm log}^{-1}|x|\chi_{
\{x:
|x|\geq 2\}}$.
The following  result generalizes estimates (\ref{toto}).

\begin{theorem}
\label{theo52}
Assume that the semigroup $e^{-tL}$ has a
kernel $p_t(x,y)$
which satisfies the upper bound (\ref{e2.1}).
If  $f\in L^{n/\alpha}({\mathbb R}^n)$ and  $L^{-\alpha/2}f<\infty $ almost
everywhere, then    $L^{-\alpha/2}f\in {{\rm BMO}_L({\mathbb R}^n)}$ with
$$
\|L^{-\alpha/2}f\|_{{\rm BMO}_L({\mathbb R}^n)}\leq
c\|f\|_{ n/\alpha}
$$
for $ 0<\alpha<n$, where  the positive constant $c$ depends only on
$\alpha$ and $n.$
\end{theorem}

Suppose  that
$T$ is a bounded operator on $L^2(\Omega).$  We say that a measurable function
$K_T\colon \Omega^2 \to {\mathbb C}$ is  the (singular) kernel  of $T$  if
     \begin{equation}\label{ker}
    \langle T f_1,f_2\rangle  =
    \int_{\Omega} T f_1(x)\overline{f_2}(x) dx  =
    \int_{\Omega}  \int_{\Omega} K_{T}(x,y) f_1(y)\overline{f_2(x)} dxdy.
    \end{equation}
for all $f_1,f_2\in C_c(\Omega)$ (for all $f_1,f_2\in C_c(\Omega)$ such that
supp $f_1 \cap $ supp $f_2 = \emptyset$ respectively).

      In order to prove Theorem~\ref{theo52}, we need the following  estimate on
the kernel  $K_{\alpha, t}(x,y)$
       of  the   operator
$({  I}-e^{-tL})L^{-\alpha/2}$ (see also Lemma 3.1 of \cite{DY3}).

\begin{lemma}
\label{lem43}  Assume that the semigroup $e^{-tL}$ has a
kernel $p_t(x,y)$
which satisfies upper bound (\ref{e2.1}).  Then for $0<\alpha<n$,
      the difference operator
$({  I}-e^{-tL})L^{-\alpha/2}$ has an associated kernel
$K_{\alpha, t}(x,y)$ which satisfies
\begin{eqnarray}
|K_{\alpha, t}(x,y)|\leq \frac{c}{ |x-y|^{n-\alpha} }\frac{t}{|x-y|^2}
\label{e3.2}
\end{eqnarray}
for some constant $c>0$
\end{lemma}
\begin{proof}
Note that
$${  I}-e^{-tL}=\int_0^t\frac{d}{dr}e^{-rL}dr=-\int_0^tLe^{-rL}dr.
$$
Hence by (\ref{e3.1})
\begin{eqnarray*}
({  I}-e^{-tL})L^{-\alpha/2}
=\frac{1}{\Gamma(\alpha/2)}\int_0^t\int_0^{\infty}
\Big(v\frac{d}{du}e^{-vL}\Big)\Big|_{v=r+s} \ \frac{1}{ r+s}
\frac{dsdr}{s^{-{\alpha/2} +1}}.
\end{eqnarray*}
By  Lemma~2.5
of \cite{CD},   the kernel of the operator $v\frac{d}{dv}e^{-vL}$
has  Gaussian upper
bound  (\ref{e2.1}).
Hence,   the   operator
$({  I}-e^{-tL})L^{-\alpha/2}$ has an associated kernel
$K_{\alpha, t}(x,y)$ which satisfies
\begin{eqnarray*}
|K_{\alpha, t}(x,y)|&\leq & c \int_0^t\int_0^{\infty}
\frac{1}{{ (r+s)^{{n}/{2}} }} e^{-c_1\frac{{|x-y|^2}}{r+s}}  \ \frac{1}{ r+s}
\frac{dsdr}{ s^{-{\alpha/2} +1}}\\
&\leq & c \int_0^t \int_0^r
\frac{1}{{ (r+s)^{{n}/{2}} }} e^{-c_1\frac{{|x-y|^2}}{ r+s}}   \ \frac{1}{r+s}
\frac{dsdr}{ s^{-{\alpha/2} +1}}\\
&&+ c \int_0^t \int_r^{\infty}
\frac{1}{{ (r+s)^{{n}/{2}} }} e^{-c_1\frac{{|x-y|^2}}{ r+s}}  \ \frac{1}{ r+s}
\frac{dsdr}{ s^{-{\alpha/2} +1}}\\
&=&{\rm   I+II}.
\end{eqnarray*}
Let us estimate  term I.   Note that  $0<s<r$. We have
\begin{eqnarray*}
      {\rm  I}
    \leq  c\int_0^t\int_0^r
      r^{-{n/2}} e^{-c_2\frac{{|x-y|^2}}{r}}
\frac{dsdr}{ rs^{-{\alpha/2} +1}}
&=&\frac{c}{ |x-y|^{n-\alpha}} \int_0^{t/|x-y|^2}
      r^{{(\alpha-n-2)/2}} e^{- {c_2  r^{-1}}}
{ dr  }\\
&\leq& \frac{c}{|x-y|^{n-\alpha}}  \frac{t}{ |x-y|^2},
\end{eqnarray*}
where the last inequality follows from   $r^{{(\alpha-n-2)/2}}
e^{-cr^{-1}}\leq c
$ for some positive constant $c$.   On the other hand,
using the condition   $0<\alpha<n$
we obtain
\begin{eqnarray*}
{\rm  II}
    \leq  c\int_0^t \int_r^{\infty}
s^{-{n/2}} e^{-c_2\frac{|x-y|^2}{ s}}
\frac{dsdr}{  s^{-{\alpha/2} +2}}
&\leq &\frac{ct}{ |x-y|^{n+2-\alpha}}
\int_0^{\infty}
      s^{{(\alpha-n-4)/2}} e^{- {c_2  s^{-1}}}
{ ds  }
\\
&\leq& \frac{c}{ |x-y|^{n-\alpha}}  \frac{t}{ |x-y|^2}.
\end{eqnarray*}
Therefore, condition (\ref{e3.2}) is satisfied and   the proof of
Lemma~\ref{lem43}
is complete.
\end{proof}
\begin{proof}[Proof of Theorem~\ref{theo52}]
      In virtue of the definition of ${\rm BMO}_L({\mathbb R}^n)$,
it suffices to   prove there exists a constant $C>0$ such that
that for any ball
$B(x,{r})$ with radius ${r}$ centered at
$x$ 
\begin{eqnarray}
\frac{1}{ |{B(x,r)}|}
\int_{{B(x,r)}}
|(I-e^{-{r}^2L})L^{-\alpha/2}f(y)|dy\leq C\| f\|_{L^{ n/\alpha}}
\label{e3.3}
\end{eqnarray}
for all $f\in L^{n/\alpha}({\mathbb R}^n)$. Set
$f_1(y)=f(y)$ if $|x-y|\le 2r$ and $f_1(y)=0$ otherwise. Next, put
$f_2=f-f_1$. Note  that
\begin{eqnarray*}
\frac{1}{ |B|}\int_{B(x,r)}
|(I-e^{-{r}^2L})L^{-\alpha/2}f(y)|dy
    &\leq& \frac{1}{ |B|}\int_{B(x,r)}
|(I-e^{-{r}^2L})L^{-\alpha/2}f_1(y)| dy
\\
&&\!\!+\frac{1}{ |B|}\int_{B(x,r)}
|(I-e^{-{r}^2L})L^{-\alpha/2}f_2(y)| dy\\
&=&\!\!{\rm  I+ II},
\end{eqnarray*}
where $|B|=|B(x,r)|$.

To estimate the first term  note that, by H\"older's
inequality
$\| f_1\|_{L^p}\leq
c|{B(x,r)}|^{ {1/ p}-{\alpha/ n}} \| f\|_{L^{n/\alpha}}.
$
for all $1<p<{n/\alpha}$.
Next,  set ${1/ q}={1/ p}-\alpha/n$.
     By Proposition~\ref{prop41}
\begin{eqnarray*}
{\rm  I}    &\leq&
\frac{1}{|B|^{1/q}}\|(I-e^{-{r}^2L})L^{-\alpha/2}f_1\|_{L^q}
\leq c  \frac{1}{|B|^{1/q}}   \|L^{-\alpha/2}f_1\|_{L^q}\\
&\leq&  c\frac{1}{|B|^{1/q}}\| f_1\|_{L^p}
\leq c \| f\|_{ L^{n/\alpha}}.
\end{eqnarray*}
To estimate the second term  note that if
$y\in {B(x,r)}$, then by Lemma~\ref{lem43}
\begin{eqnarray*}
\Big|({I}-e^{-{r}^2L})L^{-\alpha/2}  f_2 (y)\Big|
&\leq&\int_{ B(x,2r)^c} |K_{\alpha,
{r}^2}(y,z)| |f(z)|dz  \\
&\leq&c\sum_{k=1}^{\infty} \int_{ 2^k {r}\leq |x-z|< 2^{k+1} {r}} \frac{1}{
|x-z|^{n-\alpha}} \frac{{r}^2}{ |x-z|^2}
|f(z)|dz\\
&\leq&c\sum_{k=1}^{\infty}  2^{-2k} \frac{1}{ {|B(x,r2^{k+1})|}^{1-\alpha/n}}
\int_{ B(x,r 2^{k+1})}|f(z)|dz\\
&\leq&c\sum_{k=1}^{\infty}  2^{-2k}\| f\|_{ L^{n/\alpha}}
    \leq c \| f\|_{ L^{n/\alpha}}.
\end{eqnarray*}
Combining the above estimates, we obtain (\ref{e3.3}).
\end{proof}

\noindent
{\bf Remarks.}
{\rm (i)}\ Under the extra assumption that   for each
$t>0$, the kernel
$p_t(x,y)$ of
$e^{-tL}$  is  a  H\"older continuous function  in $x,$ it can be proved that
for  $f\in L^{n/\alpha}({\mathbb R}^n)$, either $
L^{-\alpha/2}f
\equiv \infty$ or $L^{-\alpha/2}f
\in {\rm BMO}_L({\mathbb R}^n)$ with
$\|L^{-\alpha/2}f\|_{{\rm BMO}_L({\mathbb R}^n)}\leq c\|f\|_{L^{n/\alpha}}.
$

We leave the details of the proof to the reader.

\medskip

      {\rm (ii)} We now give a list of examples of operators $L$ satisfying the
assumptions  in
      Proposition~\ref{prop41} and Theorem~\ref{theo52}.

\medskip
\ \ \ \ \ ($\alpha$)\ The operator ${\Delta}_N$,
${\Delta}_D$ or ${\Delta}_{DN} $
as in Section 2.3;

\medskip
\ \ \ \  ($\beta$)\
Let $V\in L^1_{\rm loc}({\mathbb R}^n)$ be a nonnegative function on
${\mathbb R}^n$ ($n\geq
3$). The Schr\"odinger operator with potential $V$ is defined  by
\begin{eqnarray}
L=-{\Delta}+V(x)\ \ \ \ {\rm on}\ {\mathbb R}^n.
\label{e6.8}
\end{eqnarray}
     From the Feynman-Kac formula, it is well-known that the    kernels
$p_t(x,y)$
of the semigroup
$e^{-tL}$  satisfy the estimate
\begin{eqnarray}
0\leq p_t(x,y)\leq \frac{1}{{ (4\pi t)^{n/2} }} e^{-\frac{|x-y|^2}{4t}}.
\label{e6.9}
\end{eqnarray}
      However,
       unless  $V$ satisfies additional conditions, the heat kernel can be
a discontinuous function of the
space variables and the H\"older continuous estimates may fail to
hold. See, for example, \cite{Da}.

We note that the
corresponding result in Theorem 1 of \cite{DGMTZ} is a special case
of  Theorem 5.2.

\medskip

\ \ \ \ ($\gamma$)\ Let $A=((a_{ij}(x))_{1\leq i,j\leq n} $ be an
$n\times n$ matrix
with  complex   entries
$a_{ij}\in L^{\infty}({\mathbb R}^n)$ satisfying $  \lambda
|\xi|^2\leq {\rm Re}\sum
a_{ij}(x)\xi_i\xi_j$ for all $x\in{\mathbb R}^n, \xi=(\xi_1, \xi_2, \cdots,
\xi_n)\in {\mathbb C}^n$ and some $\lambda>0$.
       Let $T$ be the    divergence form operator
$$
Lf \equiv -{\rm div} (A \nabla f),
$$
which we interpret in the usual weak sense via a sesquilinear form.

It is known that  Gaussian bound (\ref{e2.1}) on the heat kernel $e^{-tL}$
is  true when $A$ has real entries,
or when $n=1$, $2$ in the case of complex entries. See, for example,
\cite{AT}.

\bigskip

\subsection{Properties of fractional powers of  Neumann Laplacian
on ${\mathbb R}$}

The following example complements Theorems~\ref{theo41}~and~\ref{theo52}.
It also provides a convincing justification of introduction of the
${\rm BMO}_L$ spaces.

\begin{example}
\label{theo45}  Let ${\Delta}_N$ be  the Neumann Laplacian  on
$\RR $.  Then,
there exists  a function $f\in L^{1/\alpha}({\mathbb R})$  such that
${\Delta}_N^{-\alpha/2}f(x)<\infty$ for almost every $x\in {\mathbb R}$,
${\Delta}_N^{-\alpha/2}f\in {{\rm BMO}_{{\Delta}_N}({\mathbb R})}$
and
\begin{eqnarray}
\|{\Delta}_N^{-\alpha/2}f\|_{{\rm BMO}_{{\Delta}_N}({\mathbb
R})}\leq c\|f\|_{L^{ n/\alpha}}.
\label{e3.6}
\end{eqnarray}
However,
        $
{\Delta}_N^{-\alpha/2}f\not\in {\rm BMO}({\mathbb R}).
$
\end{example}
\begin{proof} For any $0<\alpha<1$,  we let
\begin{eqnarray}
f(x)=-\frac{1}{x^{\alpha}{\rm log} x}\chi_{\{  0<x\leq 1/2\}}(x).
\label{e3.7}
\end{eqnarray}
Then
$$
\int_{\mathbb R} |f(y)|^{1/\alpha}dy=
\int_0^{1/2}\frac{1}{ y ({\rm log}  y^{-1})^{1/\alpha}}
dy= (1-\alpha)\alpha^{-1}\big({\rm log} 2\big)^{{1/\alpha}-1} <\infty.
$$
This proves that  $f\in  L^{1/\alpha}({\mathbb R})$. It can be
verified that ${\mathcal
I}_{\alpha}f(x)<\infty$ a.e..  Also, we have that
${\Delta}_N^{-\alpha/2}f<\infty$ a.e.. Hence,

\smallskip
(a)  ${\mathcal I}_{\alpha}f\in{\rm BMO}({\mathbb R})$ with
$\|{\mathcal I}_{\alpha}f\|_{{\rm BMO}({\mathbb R})}\leq
c\|f\|_{L^{n/\alpha}}$. See
page 221 of \cite{To}.

\smallskip

(b)    By
Theorem~\ref{theo52}, we have that
${\Delta}_N^{-\alpha/2}f
\in {\rm BMO}_{{\Delta}_N} $  with  estimate (\ref{e3.6}).

\smallskip

We now prove ${\Delta}_N^{-\alpha/2}f\not\in {\rm BMO}({\mathbb R}).$
Denote by  $k^N_{\alpha}(x,y)$
      the kernel of   the fractional
powers
${\Delta}_N^{-\alpha/2}$ of ${\Delta}_N$.
By (\ref{e2.19}) and
      (\ref{e3.1})
\begin{eqnarray}
k^N_{\alpha}(x,y)=
\frac{1}{ \gamma(\alpha) }\Big (\frac{1}{|x-y|^{1-\alpha} }
+ \frac{1}{|x+y|^{1-\alpha}}\Big) H(xy),
\label{e3.8}
\end{eqnarray}
where $H$ is the Heaviside function (\ref{e2.20}).
By (\ref{e3.8})
\begin{eqnarray}
{\Delta}_N^{-\alpha/2}f(x)=
\left\{\begin{array}{ll}
0 & \ \ \ \ {\rm if}\ x\leq 0;\\[5pt]
{\mathcal I}_{\alpha}
({  f_e})(x) & \ \ \
\ {\rm if}\ x>0,
\end{array}
\right.
\label{e3.9}
\end{eqnarray}
where $f_e \in L^{1/\alpha}({\mathbb R})$ is given by the formula
$
{  f_e(x)}= -\frac{1}{|x|^{\alpha}{\rm log} |x|}
\chi_{\{    |x|\leq 1/2\}}(x).
$

For any $k\geq 5$,    we denote  $Q_k=[-{1/ k}, {1/ k}]$.
Next if
$0<x<y<1/2$, then $|x-y|<|y|$. Hence
\begin{eqnarray*}
{\Delta}_N^{-\alpha/2}f(x)&=&
\frac{1}{ \gamma(\alpha) }\int_{-{1/2}}^{{1/ 2}}  \frac{1}{|x-y|^{1-\alpha}}
{  f_e(y)}dy\\
&\geq&
-\frac{1}{ \gamma(\alpha) }\int_x^{{1/ 2}}  \frac{1}{|x-y|^{1-\alpha}}
      \frac{1}{y^{\alpha}{\rm log} y}dy\\
&\geq &
-\frac{1}{\gamma(\alpha) }\int_x^{{1/ 2}}
      \frac{1}{y {\rm log} y}dy\\
&\geq &
      \frac{1}{\gamma(\alpha) } \Big({\rm log}\ \!({\rm log} \frac{1}{ x})
-  {\rm log}\ \!({\rm log} 2)\Big),
\end{eqnarray*}
which yields
\begin{eqnarray*}
m_{Q_k}({\Delta}_N^{-\alpha/2}{f})&=&\frac{1}{|{Q_k}|}\int_{Q_k}
{\Delta}_N^{-\alpha/2}f(y)dy\\ &\geq& \frac{k}{2\gamma(\alpha)}\int_0^{1/ k}
\Big({\rm log}\ \!({\rm log} \frac{1}{y}) -  {\rm log}\ \!({\rm log}
2)\Big)dy\\
&\geq& \frac{1}{2\gamma(\alpha)}
\Big({\rm log}\ \!({\rm log} k) -  {\rm log}\ \!({\rm log} 2)\Big).
      \end{eqnarray*}
Therefore, from  (\ref{e3.9}) we obtain
\begin{eqnarray*}
\frac{1}{|{Q_k}|}\int_{{Q_k}}&&\hspace{-1cm}
|{\Delta}_N^{-\alpha/2}f(x)-m_{Q_k}({\Delta}_N^{-\alpha/2}f)|dx\\
&=&\frac{k}{2}\int_0^{{1/ k}}
|{\Delta}_N^{-\alpha/2}f(x)-m_{Q_k}({\Delta}_N^{-\alpha/2}f)|dx
+\frac{k}{2}\int_{-{1/ k}}^0 |m_{Q_k}({\Delta}_N^{-\alpha/2}f)|dx\\
&\geq &\frac{1}{2 }|m_{Q_k}({\Delta}_N^{-\alpha/2}f)|\\
&\geq& \frac{1}{4\gamma(\alpha)}
\Big({\rm log}\ \!({\rm log} k) -  {\rm log}\ \!({\rm log} 2)\Big).
\end{eqnarray*}
Note that the last term in the above inequality tends to
$\infty$  as $ k\rightarrow \infty$. Hence
$$
\sup_Q \frac{1}{ |Q|}\int_Q
|{\Delta}_N^{-\alpha/2}f(x)-m_Q({\Delta}_N^{-\alpha/2}f)|dx=\infty,
$$
where the supremum  is taken over all cubes $Q$ of ${\mathbb R}$. Therefore
${\Delta}_N^{-\alpha/2}f\not\in {\rm BMO}({\mathbb R}).$
\end{proof}

\bigskip

\noindent
{\bf Remark.}  Example~\ref{theo45} shows that for the Neumann
Laplacian ${\Delta}_N$ on the real line
${\mathbb R}$,
the ${\rm BMO}_{{\Delta}_N}({\mathbb R})$ space is
considered as a natural substitute for classical BMO space
to study the end-point boundedness of the fractional powers
${\Delta}_N^{-\alpha/2}$.

\section{${\rm BMO}_L$ estimates of imaginary powers and maximal functions.}

In this section we apply the technique of ${\rm BMO}_L$ spaces
to discuss optimal $L^p$ estimates for the imaginary powers of
the operator $L$. We refer readers to \cite{CM, G} for related results
concerning imaginary powers of self-adjoint operators.
  
Let us recall that if $L$ is a self-adjoint positive definite
operator on~$L^2({\mathbb R}^n)$. Then $L$ admits  the  spectral resolution:
$$
L= \int_0^{\infty}\lambda d E_L(\lambda),
$$
where the $E_L(\lambda)$ are spectral projectors. For any bounded
Borel function~$F\colon [0, \infty) \to {\mathbb C}$, we define the
operator $F(L)$ by the formula
\begin{equation}\label{spe}
F(L)=\int_0^{\infty}F(\lambda) d E_L(\lambda).
\end{equation}
In particular
$$
L^{is}=\int_0^{\infty} t^{is} dE(t).
$$
By spectral theory  $\|L^{is}\|_{L^2 \to L^2}=1$ for all $s\in{\mathbb R}$.
    In the following theorem we   obtain
sharp estimates for the $L^\infty \to {\rm BMO}_L$ norm of the   operators
$L^{is}$.

\begin{theorem}\label{theo61}
Assume that the heat
kernel $p_t(x,y)$ corresponding to   the self-adjoint  operator $L$ satisfies upper
bound (\ref{e2.1}).
Then
$$
\|L^{is}f\|_{{\rm BMO}_L({\mathbb R}^n)}\leq
c      (1+|s|)^{n/2}   \|f\|_{L^\infty}
$$
for all $s\in {\mathbb R}$.
\end{theorem}
\begin{proof}
It is enough to   show
that for any ball
$B(x,{r})$ with radius ${r}$ centered at
$x$, there exists a constant $C>0$ such that
\begin{eqnarray}
\frac{1}{ |{B(x,r)}|}
\int_{{B(x,r)}}
|(I-e^{-{r}^2L})L^{is}f(y)|dy\leq c (1+|s|)^{n/2}   \| f\|_{L^ \infty}.
\label{e6.3}
\end{eqnarray}
To prove (\ref{e6.3}), for any $f\in L^{\infty}({\mathbb R}^n)$,
we set
$\theta=(1+|s|)^{-1/2}$,
$f_1(y)=f(y)$ if $|x-y|\le \theta^{-1}r$ and $f_1(y)=0$ otherwise.
Next, we put $f_2=f-f_1$. Note  that
\begin{eqnarray*}
\frac{1}{ |B|}\int_{B(x,r)}
|(I-e^{-{r}^2L})L^{is}f(y)|dy
    &\leq& \frac{1}{ |B|}\int_{B(x,r)}
|(I-e^{-{r}^2L})L^{is}f_1(y)| dy
\\
&&\!\!+\frac{1}{ |B|}\int_{B(x,r)}
|(I-e^{-{r}^2L})L^{is}f_2(y)| dy\\
&=&\!\!{\rm  I+ II},
\end{eqnarray*}
where $|B|=|B(x,r)|$.
To estimate the term I we note that, by H\"older's
inequality
$$
\| f_1\|_{L^2} \leq
|{B(x,\theta^{-1}r)}|^{1/2} \| f\|_{ L^\infty}
\le |B(x,r)|^{1/2}\theta^{-n/2}\| f\|_{ L^\infty}
=|B|^{1/2}(1+|s|)^{n/2}\| f\|_{ L^\infty}.
$$
Then
\begin{eqnarray*}
{ \rm I}    \leq
|B|^{-1/2}\|(I-e^{-{r}^2L})L^{is}f_1\|_{L^2}
\leq c|B|^{-1/2}\|L^{is}f_1\|_{L^2} \\
\leq  c|B|^{-1/2}\| f_1\|_{L^2}
\leq c (1+|s|)^{n/2}\| f\|_{L^\infty}.
\end{eqnarray*}
To estimate the term  II we note that if
$y\in {B(x,r)}$, then
\begin{eqnarray*}
\Big|({I}-e^{-{r}^2L})L^{is}  f_2 (y)\Big|
&\leq&\int_{ B(x,\theta^{-1}r)^c} |K_{is,
{r}^2}(y,z)| |f(z)|dz  \\
&\leq&   \| f\|_{L^\infty}
\sup_{x\in \Omega,\, r>0}\int_{ B(x,\theta^{-1}r)^c} |K_{is,{r}^2}(x,z)|dz,
\end{eqnarray*}
where $K_{is,{r}^2}(y,z)$ is the kernel of the operator
$({I}-e^{-{r}^2L})L^{is}$.
Hence the proof of Theorem~\ref{theo61} reduces to the following
Lemma.
\end{proof}
\begin{lemma}\label{ja}
Assume that $L$ is a self-adjoint operator and its heat kernel
$p_t(x,y)$ satisfies
   the Gaussian bound  (\ref{e2.1}). Then the associated kernel
$K_{is,{r}^2}(x,z)$
of the operator $(I-e^{-{r}^2L})L^{is}$ satisfies
$$
    \int_{ B(x,\theta^{-1}r)^c} |K_{is,{r}^2}(x,z)|dz  \le c(1+|s|)^{n/2}
$$
for all $s \in {\mathbb R}$ and $r>0$.
\end{lemma}
The proof of Lemma~\ref{ja} is a minor modification of the proof of
estimates (17) of \cite{SW}. We leave the details to the reader.

Theorem~\ref{theo61} applied to the standard Laplace operator gives
the following estimates.
\begin{cor}\label{c62}
If $\Delta$ is the standard Laplace operator acting  on ${\mathbb R}^n$
then
\begin{equation}\label{lap}
\|\Delta^{is}f\|_{{\rm BMO}({\mathbb R}^n)}\leq
c      (1+|s|)^{n/2}   \|f\|_{L^\infty}
\end{equation}
for all $s\in {\mathbb R}$.
\end{cor}
\begin{proof}
Corollary~\ref{c62} is a straightforward consequence of
Theorem~\ref{theo61} and the equivalence of the classical $\rm{BMO}$ space and
${\rm BMO}_\Delta$.
\end{proof}

\medskip

\noindent
{\bf Remark.}
For the standard Laplace operator one can explicitly compute
the kernel  $|K_{is,{r}^2}(x,z)|$ and check that
    $\int_{ B(x, r)^c} |K_{is,{r}^2}(x,z)|dz  \ge c(1+|s|)^{n/2}\log(1+|s|)$.
See \cite{SW}. Hence one has to replace $B(x,2r)$ by  $  B(x,\theta^{-1}r)^c$
to obtain estimates without the additional logarithmic term.
As in \cite{SW} (Theorem 1) one can show that the norm of
    $\|\Delta^{is}\|_{L^\infty \to {\rm BMO}({\mathbb R}^n)} \ge 
c(1+|s|)^{n/2}$.
Hence the estimates in Theorem~\ref{theo61} and Corollary~\ref{c62} are
sharp. Even for the Laplace operator, our   estimate  (\ref{lap}) is
stronger than
   any
other known estimates of
$L^\infty \to {\rm BMO}$ norm of the imaginary powers of the Laplace operator.

\bigskip

Theorem~2 of \cite{SW} says that if $L$ satisfies
assumption of Theorem~\ref{theo61} then the following
estimates of the weak type $(1,1)$ norm of the imaginary powers
of $L$ holds
\begin{equation}
\label{imp}
\|L^{is}\|_{L^1 \to L^{1,\infty}} \le c (1+|s|)^{n/2}
\end{equation}

Note, however, that the week type $(1,1)$ norm is not subadditive
so despite its name is
not a norm. Whereas  $\|\cdot\|_{L^\infty \to {\rm BMO}_L}$, the norm
of linear operators form  $L^\infty$ to ${\rm BMO}_L$, is a proper
norm. This difference is crucial for the results which we discuss next.

Suppose that $F \colon \,   {\mathbb R}\to  {\mathbb C}$. Let us
recall that the Mellin transform of the function $F$ is defined by
$$
m(u)=\frac{1}{2\pi}\int_0^\infty
F(\lambda)\lambda^{-1-iu}d\lambda,\quad u\in{\mathbb R}.
$$
Moreover the inverse transform is given by the following formula
$$
F(\lambda)=\int_{\mathbb R} m(u)\lambda^{iu}du,\quad \lambda\in[0,\infty).
$$
Next we define the maximal operator $F^*(L)$ by the formula
$$
F^*(L)f(x)=\sup_{t>0}|F(tL)f(x)|,
$$
where $f\in L^p(\Omega)$ for some $1\le p \le \infty$.
\begin{cor}
Assume that $L$ is a self-adjoint operator acting on $L^2({\mathbb
R}^n)$,  and that the heat
kernel $p_t(x,y)$ of the operator $L$ satisfies upper
bound {\rm (\ref{e2.1})}.
Suppose also that $F \colon \,   {\mathbb R}\to  {\mathbb C}$ is a
bounded Borel
function    such that
$$
   \int_{\mathbb
R}|m(u)|(1+|u|)^{n/2}du  =C_{F,n} < \infty
$$
where  $m$ is the Mellin transform  of $F$. Then $F(L)$ and $F^*(L)$  are
bounded operators from $L^\infty$ to ${\rm BMO}_L$ and
$$
\|F(L)\|_{L^\infty \to {\rm BMO}_L} \le
\|F^*(L)\|_{L^\infty \to {\rm BMO}_L}
\le c C_{F,n}.
$$
\end{cor}
\begin{proof}
Note that
\begin{eqnarray*}
F(tL)&=&\int_0^{\infty}F(t\lambda) d E_L(\lambda)
=\int_0^{\infty} \int_{\mathbb R} m(u)(t\lambda)^{iu}du \,  d E_L(\lambda)\\
&=&\int_{\mathbb R}\int_0^{\infty} m(u)(t\lambda)^{iu}  d E_L(\lambda) du=
\int_{\mathbb R} m(u) t^{iu}L^{iu}du.
\end{eqnarray*}
Hence
$$
\sup_{t>0}|F(tL)f(x)| \le  \int_{\mathbb R} |m(u)|   |L^{iu}f(x)|  du
$$
and
$$
\|F^*(L)f  \|_{{\rm BMO}_L}
    \le \int_{\mathbb R} |m(u)|    \|f\|_{L^\infty}
\|L^{iu}\|_{L^\infty \to {\rm BMO}_L}du
\le c \|f\|_{L^\infty} \int_{\mathbb R}|m(u)|(1+|u|)^{n/2}du.
$$
The inequality $\|F(L)\|_{L^\infty \to {\rm BMO}_L} \le
\|F^*(L)\|_{L^\infty \to {\rm BMO}_L}$ is an obvious consequence of
the definition of $F^*(L)$.
\end{proof}

\section{${\rm BMO}_L$ estimates for spectral multipliers
of self-adjoint operators}

In this section we discuss an application of ${\rm BMO}_L(\Omega)$ technique
to the theory of  H{\"o}rmander spectral multipliers.
In the sequel if $F(L)$ is the operator defined by (\ref{spe}) then
by $K_{F(L)}$ we denote
the kernel associated with   $F(L)$. See (\ref{ker})
of \cite{DY1}.

\begin{theorem} \label{CZ}
Suppose that $\| F \|_{L^\infty} \le C_1$, and that
\begin{equation}\label{lap1}
\sup_{r>0}\sup_{y\in \Omega}\int_{B(y,r)^c}
|K_{F(L)(I-e^{-r^2L})}(x,y)| dx \leq C_1.
\end{equation}
Then
$$
    \| F(L)\|_{L^\infty \to {\rm BMO}_{L}} \le c C_1.
$$
\end{theorem}
\begin{proof}
We note again that   it is enough to   show
that for any ball
$B(x,{r})$ with radius ${r}$ centered at
$x$, there exists a constant $C>0$ such that
\begin{eqnarray}
\frac{1}{ |{B(x,r)}|}
\int_{{B(x,r)}}
|(I-e^{-{r}^2L})F(L)f(y)|dy\leq c C_1  \| f\|_{L^ \infty}.
\label{e6.3n}
\end{eqnarray}
To prove (\ref{e6.3n}) for any $f\in L^{\infty}({\mathbb R}^n)$ we set
$f_1(y)=f(y)$ if $|x-y|\le 2r$ and $f_1(y)=0$ otherwise.
Next, we put $f_2=f-f_1$. Note  that
\begin{eqnarray*}
\frac{1}{ |B|}\int_{B(x,r)}
|(I-e^{-{r}^2L})F(L)f(y)|dy
    &\leq& \frac{1}{ |B|}\int_{B(x,r)}
|(I-e^{-{r}^2L})F(L)f_1(y)| dy
\\
&&\!\!+\frac{1}{ |B|}\int_{B(x,r)}
|(I-e^{-{r}^2L})F(L)f_2(y)| dy\\
&=&\!\!{\rm  I+ II},
\end{eqnarray*}
where $|B|=|B(x,r)|$.
To estimate the term I we note that, by H\"older's
inequality
$$
\| f_1\|_{L^2} \leq
|{B(x,2r)}|^{1/2} \| f\|_{ L^\infty}
\le   c   |{B(x,2r)}|^{1/2} \| f\|_{ L^\infty} .
$$
Then
\begin{eqnarray*}
{ \rm I}    \leq
|B|^{-1/2}\|(I-e^{-{r}^2L})F(L)f_1\|_{L^2}
\leq c|B|^{-1/2}\|F(L)f_1\|_{L^2} \\
\leq  c|B|^{-1/2}C_1 \| f_1\|_{L^2}
\leq c C_1\| f\|_{L^\infty}.
\end{eqnarray*}
To estimate the term II we note that if
$y\in {B(x,r)}$, then
\begin{eqnarray*}
\Big|({I}-e^{-{r}^2L})F(L) f_2 (y)\Big|
&\leq&\int_{ B(y, r)^c} |K_{( {I}-e^{-{r}^2L})F(L) }(y,z)| |f(z)|dz  \\
&\leq&   \| f\|_{L^\infty}
\sup_{x\in \Omega,\, r>0}\int_{ B(y,r)^c} |K_{  (
{I}-e^{-{r}^2L})F(L)}(y,z)|dz \\
&\le& cC_1 \|
f\|_{L^\infty}
\end{eqnarray*}
\end{proof}

In the standard theory of H{\"o}rmander spectral multipliers one usually
begins with proving weak type
$(1,1)$ estimates for a spectral multiplier $F(L)$. Next $F(L)$ is
bounded on $L^2$ by
the spectral theorem so continuity of the operator $F(L)$ on $L^p$
spaces for $1< p < \infty$ follows from the Marcinkiewicz
interpolation theorem.
One can use Theorem~\ref{CZ} and  Proposition~\ref{prop3} to obtain
an alternative proof of boundedness of $F(L)$ on an $L^p$ space for
$1< p < \infty$.
Of course continuity of $F(L)$ as an  operator from $L^\infty$ to ${\rm BMO}_L$
is   of independent interest even if we already know that $F(L)$ is
of weak type
$(1,1)$.

   The  H{\"o}rmander type spectral multipliers is a very broad subject.
For example such multipliers were studied in \cite{Ale, Ch, DOS, He, MM, MS}.
One can use Theorem~\ref{CZ} to show that  all spectral multipliers
of weak type $(1,1)$ which are discussed in  \cite{Ale, Ch, DOS, He, MM, MS}
are also bounded from $L^\infty$ to ${\rm BMO}_L$.
As an example we discuss the following ${\rm BMO}_L$ versions of
Theorem~3.1 of \cite{DOS}. Let us recall that if $F \colon {\mathbb R} \to
{\mathbb C}$ then $\|F\|_{W^{p}_s}=\|(I+\Delta)^{n/2}F\|_{L^p({\mathbb R})}$.

\begin{theorem}\label{des}
Suppose that $L$ is a self-adjoint operator  acting on $L^2(\Omega)$,
$\Omega \subset
    {\mathbb R}^n$
and that the heat kernel $p_t(x,y)$ of
   $L$ satisfies the Gaussian bound  {\rm (\ref{e2.1})} and that
$\eta \in C^{\infty}_c({\mathbb R}_+)$. Then  for every  $s>n/2$  and
for all Borel bounded function
     $F$ such that  $\sup\limits_{t> 0} \Vert \eta \, \delta_t F
\Vert_{W^\infty_s}<\infty$ the operator $F(L)$ is bounded on $L^p(\Omega)$
for all $1<p<\infty$. Moreover
    \begin{equation}\label{hor:m1}
       \|F({L})\|_{L^{\infty}   \to {\rm BMO}_L}
       \le  C_s  \Big(\sup_{t >   0} \Vert \eta
       \, \delta_t F \Vert_{W^\infty_s}\Big)
       \end{equation}
for all $s>n/2$.
\end{theorem}
\begin{proof}
Note that by  \cite{DOS} ((4.19) and   Remark~1), we have
$$
\sup_{r>0}\sup_{y\in \Omega}\int_{B(y,r)^c}
|K_{F(L)(I-e^{-r^2L})}(x,y)| dx \leq  C_s  \Big(\sup_{t >   0} \Vert \eta
       \, \delta_t F \Vert_{W^\infty_s}\Big).
$$
Hence Theorem~\ref{des} is a straightforward consequence of Theorem~\ref{CZ}.
\end{proof}

{\bf Acknowledgments.}  The authors  thank   A. McIntosh  for
      helpful suggestions and discussions.

\bigskip

\bigskip

{\footnotesize

DEPARTMENT OF MATHEMATICS, ZHONGSHAN  UNIVERSITY, GUANGZHOU, 510275, P.R. CHINA

{\it E-mail address}: stsddg@mail.sysu.edu.cn

\bigskip

{ DEPARTMENT OF MATHEMATICS, MACQUARIE UNIVERSITY, NSW 2109, AUSTRALIA}

{\it E-mail address}: duong@ics.mq.edu.au

\bigskip

DEPARTMENT OF MATHEMATICAL SCIENCES, NEW MEXICO STATE UNIVERITY, P.O.
BOX 30001, LAS
CRUCES, NM, 88003-8001

{\it E-mail address}: asikora@nmsu.edu

\bigskip

DEPARTMENT OF MATHEMATICS, MACQUARIE UNIVERSITY, NSW 2109, AUSTRALIA {  AND}
DEPARTMENT OF MATHEMATICS, ZHONGSHAN  UNIVERSITY, GUANGZHOU, 510275, P.R. CHINA

{\it E-mail address}:  mcsylx@mail.sysu.edu.cn

}

\end{document}